\documentclass[aos,MSNbibl,seceqn,dvips]{arximspdf}
\usepackage{graphicx}
\doi{10.1214/13-AOS1126} %kopijuoti is PTS
\volume{41}
\issue{3}
\pubyear{2013}
\firstpage{1542}
\lastpage{1568}

\makeatletter

\newproclaim{remark}{Remark}[section]
\newtheorem{lemma}{Lemma}[section]
\newtheorem{corollary}{Corollary}[section]
\newtheorem{proposition}{Proposition}[section]
\newtheorem{theorem}{Theorem}[section]
\newproclaim{assumption}{Assumption}[section]

\newcommand{\cal}{\mathcal}

\makeatother

\begin{document}
\begin{frontmatter}

\title{Nonparametric regression with the scale depending on auxiliary variable}
\runtitle{Nonparametric regression}

\begin{aug}
\author[A]{\fnms{Sam} \snm{Efromovich}\corref{}\thanksref{T1}\ead[label=e1]{efrom@utdallas.edu}}
\runauthor{S. Efromovich}
\affiliation{University of Texas at Dallas}
\address[A]{Department of Mathematical Sciences\\
University of Texas at Dallas\\
Richardson, Texas 7580\\
USA\\
\printead{e1}} %adresu isvedimo komanda gale!
\end{aug}

\thankstext{T1}{Supported in part by NSF Grant DMS-09-0679 and NSA Grant H982301310212.}

% HISTORY:
\received{\smonth{9} \syear{2012}}
\revised{\smonth{4} \syear{2013}}

% ABSTRACT
%
\begin{abstract}
The paper is devoted to the problem of estimation of a univariate
component in a heteroscedastic nonparametric multiple regression
under the mean integrated squared error (MISE) criteria. The aim is to
understand how the scale function should be used for estimation of the
univariate component. It is known that the scale function does not
affect the rate of the MISE convergence, and as a result sharp
constants are explored. The paper begins with developing a
sharp-minimax theory for a pivotal model $Y=f(X) + \sigma(X,{\mathbf Z})
\varepsilon$, where $\varepsilon$ is standard normal and independent of the
predictor $X$ and the auxiliary vector-covariate ${\mathbf Z}$. It is
shown that if the scale $\sigma(x,{\mathbf z})$ depends on the auxiliary
variable, then a special estimator, which uses the scale (or its
estimate), is asymptotically sharp minimax and adaptive to unknown
smoothness of $f(x)$. This is an interesting conclusion because if the
scale does not depend on the auxiliary covariate ${\mathbf Z}$, then
ignoring the heteroscedasticity can yield a sharp minimax estimation.
The pivotal model serves as a natural benchmark for a general additive
model $Y = f(X) + g({\mathbf Z}) + \sigma(X,{\mathbf Z})\varepsilon$, where
$\varepsilon$ may depend on $(X,{\mathbf Z})$ and have only a finite fourth
moment. It is shown that for this model a data-driven estimator can
perform as well as for the benchmark. Furthermore, the estimator,
suggested for continuous responses, can be also used for the case of
discrete responses. Bernoulli and Poisson regressions, that are
inherently heteroscedastic, are particular considered examples for
which sharp minimax lower bounds are obtained as well. A numerical
study shows that the asymptotic theory sheds light on small samples.
\end{abstract}

% KEYWORDS

\end{frontmatter}

%s1 #&#
\section{Introduction}\label{sec1}

We begin the \hyperref[sec1]{Introduction} with a simple model which
will allow us to explain the setting and the problem, then formulate
studied extensions and finish with terminology used in the paper.

%s1.1 #&#
\subsection{Pivotal regression model}\label{sec1.1}

In order to set the stage for a variety of considered problems, it is
convenient to begin with a simple nonparametric regression model
%
%e1.1 #&#
%
\begin{equation}\label{equ1.1}
Y = f(X) + \sigma(X, {\mathbf Z})\varepsilon,
\end{equation}
which will serve as a pivot for all other models. In (\ref{equ1.1})
$Y$ is the response, $X$ is the univariate random predictor of interest and
${\mathbf Z}:= (Z_1,\ldots,Z_{D})$
is the vector of random auxiliary covariates,
$\sigma(x,{\mathbf z})$ is the scale function
[$\sigma^2(x,{\mathbf z})$ is also called the variance or volatility] and
$\varepsilon$ is a standard normal error independent of
$(X,{\mathbf Z})$. It is assumed that $(X,{\mathbf Z})$
has a joint density $p(x,{\mathbf z})$
supported on $[0,1]^{1+D}$, and in what follows $p(x)$ denotes the
(marginal) density of $X$.
The problem is to estimate the nonparametric regression
function $f(x)$ based on a sample of size $n$ from $(X,{\mathbf Z}, Y)$.

\begin{figure}

\includegraphics{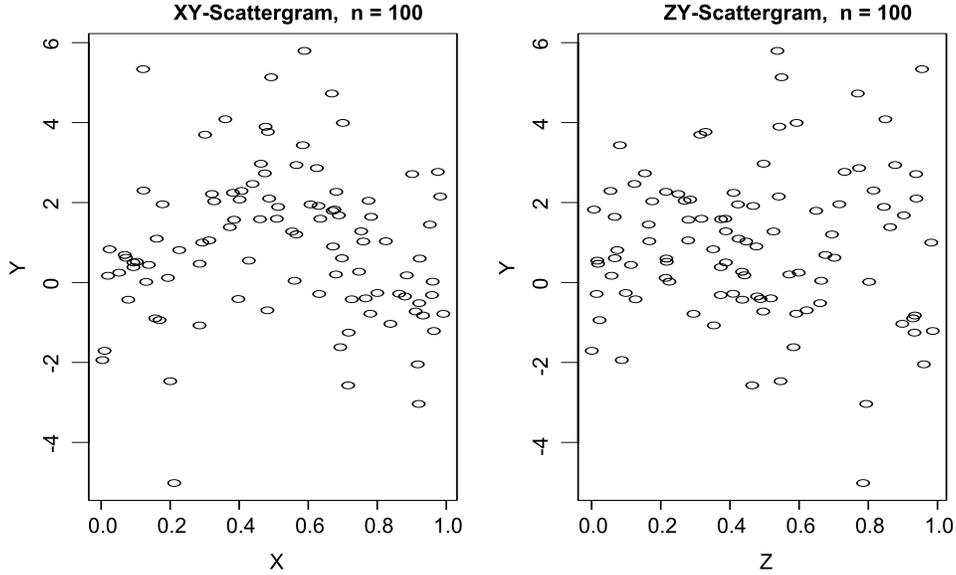}

\caption{Scattergrams for a data simulated according to
model (\protect\ref{equ1.1}) with
$D=1$.}\label{fig1}
\end{figure}

Figure~\ref{fig1}
illustrates model (\ref{equ1.1}) for a particular case $D=1$
and $n=100$ (more details will be revealed shortly). The data is volatile
(compare with ``typical'' data studied in~\cite{8,12,19,34}),
and it is difficult to visualize an underlying regression. The $XY$-scattergram
suggests a number of possible outliers, but here we do know
that these are not outliers, and they are due to heteroscedasticity
that can be observed in the $ZY$-scattergram.
Typically, for such a data with two covariates one would definitely
%
%f1 #&#
%
attempt to use a multiple or additive
regression to explain or reduce the volatility in $XY$-scattergram and
to improve
visualization of the underlying regression. However, here
we do know that there is no additive component in $z$. The only
hope to help a nonparametric estimator is to use a known (or estimated)
scale function.
But is this worthwhile to do, and if the answer is ``yes,''
then how one should proceed?
%This paper is devoted to exploring the usefulness of the scale
%in nonparametric regression estimation, we will return to discussion
%of this
%simulation in Section 4, and
Before presenting the answer, let us return to describing the studied
setting and known results.

%s1.2 #&#
\subsection{Pivotal problem}\label{sec1.2}

To be specific about smoothness of $f(x)$ and because we are going to
study minimax constants, let us assume that $f(x)$ belongs to a Sobolev
class ${\cal S}(\alpha,Q):= \{f(x)\dvtx f(x)=\sum_{j=0}^\infty
\theta_j \varphi_j(x), \varphi_0(x):=1,\break \varphi_j(x):= 2^{1/2}
\cos(\pi j x), j\geq1$, $\sum_{j=0}^\infty[1+(\pi j)^{2\alpha}]
\theta_j^2 \leq Q < \infty, x \in[0,1], \alpha\geq1\}$.
Furthermore, the risk of an estimate $\check f(x)$ is defined by the
mean integrated squared error (MISE) $E\{\int_0^1 (\check f(x) -
f(x))^2 \,dx\}$.

The above-presented discussion of a simulation exhibited in Figure~\ref{fig1}
raises the following question.
Suppose that, apart of $f(x)$, the statistician knows everything about
regression (\ref{equ1.1}). Should one use the scale function
(and correspondingly the auxiliary variable) in a regression
estimator? To warm up the reader, let
us consider several arguments against and for using the scale.
Against: (a1) A majority of nonparametric research
is devoted to rates of the MISE
convergence. For the considered setting the rate is
$n^{-2\alpha/(2\alpha+ 1)}$, and then practically all known
estimators can attain this rate without using the scale;
see~\cite{12,13,14,20}. (a2)
There is a widely held opinion that regression estimation is
``\ldots relatively insensitive to heteroscedasticity\ldots'' as discussed
in~\cite{34}. (a3) This is probably the strongest argument against
using/estimating the scale.
Let us consider a particular case $\sigma(x,{\mathbf z}) =
\sigma(x)$ and assume that $p(x)$ and $\sigma(x)$ are positive and have bounded
derivatives on $[0,1]$. Then in~\cite{7} the
following sharp minimax lower bound is established:
%
%e1.2 #&#
%
\begin{eqnarray}\label{equ1.2}
&&\inf_{\check f^*} \sup_{f \in{\cal S}(\alpha,Q)} E\biggl\{\int
_0^1 \bigl(\check f^*(x; p,\sigma,\alpha,Q)-
f(x)\bigr)^2\,dx\biggr\}
\nonumber\\[-8pt]\\[-8pt]
&&\qquad\geq P(\alpha,Q) \bigl[d_1(p,\sigma) n^{-1}
\bigr]^{2\alpha/(2\alpha+ 1)}\bigl(1+o_n(1)\bigr),\nonumber
\end{eqnarray}
where the infimum is taken over all possible $\check f^*$
based on a sample
$\{(X_1,Y_1),\break\ldots,(X_n,Y_n)\}$,
the design density $p(x)$, the scale function $\sigma(x)$ and
parameters $(\alpha, Q)$ that define the underlying Sobolev class. In
(\ref{equ1.2})
%
%e1.3 #&#
%
\begin{equation}\label{equ1.3}
P(\alpha,Q):= \bigl[\alpha/\pi(\alpha+1)\bigr]^{2\alpha/(2\alpha+ 1)}
\bigl[Q(2
\alpha+1)\bigr]^{1/(2\alpha+1)}
\end{equation}
is the Pinsker constant~\cite{30}, and
%
%e1.4 #&#
%
\begin{equation}\label{equ1.4}
d_1:= d_1(p,\sigma):= \int_0^1
\frac{\sigma^2(x)}{p(x)}\,dx
\end{equation}
is the coefficient of difficulty which is equal to one in the classical
case of the unit scale and uniform design, and
here and in what follows $o_n(1)$'s denote
generic sequences which vanish as $n \to\infty$.
Furthermore, in~\cite{7} (see also~\cite{8}) it is shown that
there exists an estimator based solely on data
(in what follows referred to as $E$-estimator) that does not estimate
the scale
$\sigma(x)$, ``ignores'' the heteroscedasticity and nonetheless attains
the lower bound (\ref{equ1.2}).
In other words, the
``ignore-heteroscedasticity'' methodology may yield
a sharp-minimax estimation.
Furthermore, according to~\cite{7,8}
the $E$-estimator performs well for small samples.

Typical arguments in favor of using/estimating the scale are as follows:
(f1) Scale affects the constant of the MISE convergence, and constants
may be more important for small samples than rates~\cite{8,27,28};
(f2) Weighted regression (with weights depending on the scale)
is a familiar remedy for heteroscedasticity \cite
{12,14,15,18,20,31,34}; (f3)
It is reasonable to believe that using the scale
may improve an estimator.

Because there are many rate-optimal estimators, to answer the
raised pivotal question it is natural
to explore a sharp-minimax estimation, that is,
estimation with best constant and rate of the MISE convergence.
It will be shown shortly that for the
model (\ref{equ1.1}) the lower bound (\ref{equ1.2}) [with the infimum
taken over all possible
$\check f^*$ based on a sample of size $n$ from $(X,{\mathbf Z},Y) $,
all nuisance functions defining the model (\ref{equ1.6}) and parameters
$(\alpha,Q)$]
still holds with $d_1$ being replaced by
%
%e1.5 #&#
%
\begin{equation}\label{equ1.5}
d:= d(p,\sigma):= \int_0^1 \frac{dx} {p(x) E\{\sigma^{-2}(X,{\mathbf
Z})|X=x\}}.
\end{equation}
The $E$-estimator, if it is na\"ively used for model (\ref{equ1.1}),
is consistent and even rate minimax, and
supremum (over the Sobolev class) of its MISE is equal to the right-hand
side of (\ref{equ1.2}) with $d_1$ being replaced by
$d_2:= E\{\sigma^2(X,{\mathbf Z})p^{-2}(X)\}$. The latter, according to
the Cauchy--Schwarz inequality,
is larger than $d$ whenever the scale depends on the auxiliary variable.

We conclude that for the scale depending on an auxiliary variable,
the $E$-estimator, which ignores heteroscedasticity, is no
longer sharp minimax. As a result, it is reasonable to
explore a regression estimator that uses the scale to attain the
sharp minimaxity. The underlying idea of the proposed estimator,
based on the developed asymptotic theory, is to use weighted responses
$w_lY_l$ with weights
\[
w_l(p,\sigma):= %\frac{\sigma^{-2}(X_l,{\mathbf Z}_l)}{
p^{-1}
(X_l) \frac{\sigma^{-2}(X_l,{\mathbf Z}_l)}{E\{ \sigma
^{-2}(X_l,{\mathbf Z}_l)|X_l\}}.
\]
Note that: $p^{-1}(X_l)$ is a well-known weight in a univariate
sharp-minimax regression~\cite{8}; If $\sigma(x,{\mathbf z}) = \sigma
(x)$, then the
weight does not depend on the scale; Given $X_l=x_l$, conditional expectation
$E\{\sigma^{-2}(x_l,{\mathbf Z}_l)|X_l=x_l\}$
is the best estimate (predictor) of $\sigma^{-2} (x_l,{\mathbf Z}_l)$ under
the MSE criteria, and the better the
estimation is, the closer the weight will be to $p^{-1}(X_l)$; In the
light of the foregoing, the
proposed weight may be of a special benefit to the case of
independent $X$ and ${\mathbf Z}$; The weights should help in dealing with
``outliers'' created by
heteroscedasticity in auxiliary covariates. To shed additional light on
the made comments,
let us return to Figure~\ref{fig1}.
The underlying model is defined in Section~\ref{sec4} where it is revealed
that the used scale is $\sigma(x,z) = \sigma(z)$ and $X$ and $Z$ are
independent.
[The interested reader can also look at the identical left diagram in
Figure~\ref{fig2} where the solid
line shows the underlying regression $f(x)$.] We can now realize that
``outliers''
in the $XY$-scattergram are created
by the heteroscedasticity in $z$ and the independence of $Z$ from $X$
which creates a
chaotic placement of ``outliers'' in the scattergram.
%Furthermore, it is intuitively clear that the proposed weights may
%help to attenuate the effect of heteroscedasticity in the auxiliary
%covariate.

%On the other hand, let us for a moment change the setting
%and, using the same Figure 1, consider a case where
%regression of $Y$ on $Z$ is of interest. The
%scattergram of interest is shown in the right diagram and it exhibits
%a pronounced heteroscedasticity. Nonetheless,
%because there is no heteroscedasticity in the ``auxiliary''
%$XY$-scattegram (and as a result no otliers in the right diagram),
%we do not need to use the scale in weighted responses for sharp
%minimax estimation. Of course,
%in this example we have an additive nuisance component, and we are
%considering this issue below.

%s1.3 #&#
\subsection{Extensions}\label{sec1.3}

The following extensions of the model (\ref{equ1.1}) will be considered:

(i) Model (\ref{equ1.1}) is a natural benchmark for a general additive model
%
%e1.6 #&#
%
\begin{equation}\label{equ1.6}
Y = f(X) + g({\mathbf Z}) + \sigma(X,{\mathbf Z}) \varepsilon,
\end{equation}
where $g({\mathbf z})$ is a nuisance $D$-dimensional additive
component integrated to zero on $[0,1]^D$.
There is a vast literature devoted to univariate additive models
\cite{14,15,17,18,21,22,23,24,33,35}, with the most advanced
sharp-minimax result
due to Horowitz, Klemela and Mammen~\cite{22}
where, for the case of a known $\sigma(x, {\mathbf z}) = \sigma$,
$g({\mathbf z}) = g_1(z_1)+\cdots+g_D(z_D)$
with differentiable univariate additive components,
and known parameters $\alpha$, $Q$ and $\sigma$,
a shrinkage estimator $\check f(x,\alpha,Q, \sigma)$ is proposed
such that for any $C > 0$,
\begin{eqnarray*}
&&
\sup_{f \in{\cal S}(\alpha, Q)}\operatorname{Pr} \biggl( (n/
d_1)^{2\alpha
/(2\alpha+ 1)}
P^{-1}(\alpha,Q)
\\
&&\hspace*{31.5pt}\qquad{}
\times E\biggl\{\int_0^1 \bigl(\check f(x,
\alpha,Q,\sigma) - f(x)\bigr)^2 \,dx| (X_1,{\mathbf
Z}_1),\ldots,(X_n,{\mathbf Z}_n)\biggr\} \\
&&\hspace*{269.5pt}\qquad> 1+C \biggr)\\
&&\quad= o_n(1).
\end{eqnarray*}
We will show shortly that
without any assumption on the structure of unknown $g({\mathbf z})$ there
exists a data-driven sharp-minimax estimator.
In other words, the presence of
a nuisance additive component $g({\mathbf z})$ affects neither minimax rate,
nor the sharp minimax constant, nor the ability of adaptive estimation.

\mbox{}\hphantom{i}(ii) It is of interest to relax the assumption
about independence between the regression error and covariates
as well as the assumption about normal distribution of the error.
It will be shown shortly that the MISE of the proposed regression
estimator still attains the
minimax lower bound\vadjust{\goodbreak} (\ref{equ1.2}), with $d_1$ being replaced by $d$, whenever
the regression error satisfies
%
%e1.7 #&#
%
\begin{equation}\label{equ1.7}
E\{\varepsilon|X,{\mathbf Z}\} = 0,\qquad E\bigl\{\varepsilon
^2|X,{\mathbf
Z}\bigr\} =
1,\qquad E\bigl\{\varepsilon^4|X,{\mathbf Z}\bigr\} < C <
\infty\qquad\mbox{a.s.}\hspace*{-28pt}
\end{equation}
To compare with a known assumption for a univariate regression,
in~\cite{7} for model (\ref{equ1.1}) with
$\sigma(x,{\mathbf z})=\sigma(x)$ the proposed adaptive estimation
assumes independence of the predictor and regression error $\varepsilon
$ plus
a finite eighth moment of the regression error.

(iii) Extension (ii) is a natural bridge to other classical
heteroscedastic models as well as to discrete responses. In this paper
Bernoulli and Poisson regressions, that are inherently heteroscedastic,
are considered. Note that these regressions create a new
issue of satisfying bona fide properties of the regression function, and
the following extension is instrumental in solving the issue.

\mbox{}\hspace*{1pt}(iv) As we shall see shortly, it is worthwhile to
replace a single Sobolev class ${\cal S}(\alpha,Q)$ by a family ${\cal
F}$ of function classes that includes Sobolev, local Sobolev
(introduced in Golubev~\cite{16}) and shrinking
(toward a pivotal regression function) Sobolev
classes as particular cases. Namely, set
%
%e1.8 #&#
%
\begin{eqnarray}\label{equ1.8}
{\cal F}&:=& {\cal F}(f_0,\rho_n,M_n,
\alpha,Q)
\nonumber\\
&:=& \Biggl\{f(x)\dvtx f(x)=\sum_{j=0}^{M_n-1} \int
_0^1f_0(u) \varphi_j(u)
\,du \varphi_j(x) I(M_n > 0) + \sum
_{j \geq M_n} \theta_j \varphi_j(x),\nonumber\\
&&\hspace*{8pt} x\in[0,1],\nonumber\\[-8pt]\\[-8pt]
&&\hspace*{11pt}\sup_{x \in[0,1]}\bigl|f_0(x)\bigr| < \infty, \int
_0^1f_0^2(x)\,dx <
\infty, \theta_j:=\int_0^1 f(u)
\varphi_j(u) \,du,
\nonumber\\
&&\hspace*{11pt}\sum_{j \geq M_n}\bigl[1+(\pi j)^{2\alpha}\bigr]
\theta_j^2 \leq Q < \infty, \sup_{x \in[0,1]}\biggl|
\sum_{j
\geq M_n} \theta_j \varphi(x)\biggr| <
\rho_n,\nonumber
\\
&&\hspace*{70pt}\alpha\geq1, 0 \leq M_n < n^{1/(2\alpha+1)}/\ln^2(n),
\rho_n > n^{-1/(2\alpha+
1)}\ln(n) \Biggr\}.\nonumber
\end{eqnarray}
Here $f_0(x)$ is a bona fide (e.g., positive for Poisson regression)
regression function which will be referred to as a \textit{pivot},
$I(\cdot)$ is the indicator and the last line in (\ref{equ1.8})
specifies restrictions on $\alpha$ and numerical sequences $\rho_n$ and~$M_n$.

%s1.4 #&#
\subsection{Comments on the family ${\cal F}$ and minimax approach}\label{sec1.4}

(a) With respect to a classical Sobolev
class ${\cal S}(\alpha, Q)$, we have
${\cal S}(\alpha,Q) = {\cal F}(0,\infty,0,\alpha,Q)$,
and if the pivot is constant $f_0(x) = C$, $C < Q^{1/2}$, then ${\cal
F}(C,\rho_n,1,\alpha,Q-C^2) \subset{\cal S}(\alpha,Q)$. As a result,
the classical Sobolev class is a particular
(not changing with $n$) member of the family.
A function $f$ from the family is not farther
than $\rho_n$ in $L_\infty$-norm from the pivot. Furthermore,
if $M_{n} > 0$, then on $M_{n}$ low frequencies the regression function
$f$
is equal to the pivot, and on higher frequencies it is not farther
than $\rho_{n}$
in $L_\infty$-norm and not farther than $([1+(\pi
M_n)^{2\alpha}]^{-1}Q)^{1/2}$ in $L_2$-norm. As a result, if either
$\rho_n$ or $M_n^{-1}$ vanishes as $n \to\infty$, the set of considered
regression functions shrinks toward the pivot. This allows us to
conclude that the family ${\cal F}$
includes local Sobolev classes shrinking in $L_2$-norm, or $L_\infty
$-norm, or in both norms
to the pivot. Two other shrinking properties are
${\cal F}(f_0,\rho,M,\alpha,Q) \subset{\cal
F}(f_0,\rho+\gamma,\break M,\alpha,Q)$ and ${\cal
F}(0,\rho,M+\gamma,\alpha,Q) \subset{\cal
F}(0,\rho,M,\alpha,Q)$, $\gamma> 0$. Let us also note that
a local Sobolev class, proposed in Golubev~\cite{16},
can be written as $f_0 + {\cal F}(0,\rho,0,\alpha,Q)$
where $f_0 \in{\cal S}(\alpha',Q')$,
$\alpha' > \alpha$ and $\rho> 0$. Furthermore,
let us note that $n^{1/(2\alpha+1)}$ is the
classical number of Fourier coefficients that should be estimated by a
rate-minimax estimator; this sheds light on the upper bound in
the last line of (\ref{equ1.8}) for considered
$M_n$. The lower bound for considered $\rho_n$ is due to a specific
least favorable prior distribution of parameters which is used in
establishing the minimax lower bound.

(b) It may be convenient to think about both the function family (\ref
{equ1.8}) and
the minimax approach in terms of the game theory.
There are three players in a minimax game: the dealer, nature and
the statistician. The game is defined by: (i) an underlying \textit
{model} [here
a regression model (\ref{equ1.6})]; (ii)~assumptions
about \textit{nuisance functions} [here the additive component
$g({\mathbf z})$,
scale $\sigma(x,{\mathbf z})$, distribution of the error $\varepsilon$
and the design
density $p(x,{\mathbf z})$]; (iii)~\textit{parameters} of a family
${\cal F}$
which defines a class of estimated functions $f(x)$ [here
$ {\cal F}$ is defined in (\ref{equ1.8}) and the
parameters are the pivotal regression $f_0(x)$, sequences $M_n$ and
$\rho_n$
and Sobolev's $\alpha$ and $Q$].
The game begins with the dealer dealing nuisance functions and
parameters of ${\cal F}$ to nature. This deal must satisfy
assumptions of the game. Then for each $n$ nature
chooses a regression function $f$ from the dealt ${\cal F}$
and generates a sample of size $n$ using $f$ and the dealt model.
The dealer and the statistician, using the sample, estimate $f$.
The dealer knows everything apart of estimated $f$, the statistician
knows the sample, all assumptions of the game plus may know some
nuisance functions (like the design density in controlled regressions or
the distribution of error in special regression models like Poisson).
Nature tries to select most difficult regression function $f$ for estimation,
and the dealer and the statistician try
to estimate it with the smallest MISE. The dealer has an advantage of knowing
the dealt ${\cal F}$ and nuisance functions, and therefore the dealer's
MISE may serve as a lower bound (benchmark) for the statistician.
%Another remark is about terminology: we may say that
%(\ref{equ1.8}) describes a particular deal,
%and that made assumptions define the totality (family) of possible
%deals.

(c) Using family (\ref{equ1.8}) of function classes in place of
a single Sobolev class allows us to answer (at least partially)
to a familiar criticism of a minimax approach
that the statistician cares only about the worst case scenario
regression from ${\cal S}(\alpha,Q)$
which can be far from an underlying regression function.
This is where introducing a pivot whose smoothness is not restricted,
together with the possibility to consider shrinking function classes, shines.

(d) The smaller a function
class is, the smaller the minimax MISE (for the dealer and the statistician)
may be.
This is where the imposed restriction [see the last line in (\ref{equ1.8})]
on the dealer's choice of deals
comes into the play. As we shall see shortly, none of the legitimate deals
(which may imply local and/or shrinking function classes)
changes a sharp lower bound known for
a classical Sobolev class ${\cal S}(\alpha,Q)$. On the other hand,
not all estimates, that are sharp minimax for Sobolev classes, are even
rate minimax for particular deals. For instance,
classical estimates based on the Pinsker smoothing,
used for a univariate
regression model in Efromovich~\cite{8} and an additive regression model
in Horowitz, Klemela and Mammen~\cite{22}, are sharp minimax for a
Sobolev class,
but not even rate minimax for ${\cal F}$ whenever pivot $f_0$ and
sequence $M_n$ are such that
$\sum_{j=1}^{M_n} j^{2\alpha}[\int_0^1 f_0(x)
\varphi_j(x) \,dx]^2 \to\infty$ as $n \to\infty$.
In other words, if the pivot is not
a Sobolev function of order $\alpha$, then the famous Pinsker smoothing
is no longer even rate minimax. We will prove this assertion in
the Appendix (see~\cite{Efr}).

%s1.5 #&#
\subsection{Terminology}\label{sec1.5}
The aforementioned approach [Section~\ref{sec1.4}(b)] allows us to
introduce the following terminology. \textit{Estimator} is a statistic
based on a sample, made assumptions and, if known, on nuisance
functions defining model (\ref{equ1.6}).
% with particular examples being $g({\mathbf z})$, $p(x,{\mathbf z})$ or
%$\sigma(x,{\mathbf z})$. For instance, in many applications it is
%assumed
%that a regression function
%the joint design density $p(x,{\mathbf z})$ is always known in
%regressions
%with controlled design,
%and $g({\mathbf z})=0$ in regression model (\ref{equ1.1}).
In what follows we will explicitly state what nuisance functions,
if any, are known. \textit{Dealer-estimator}
knows everything about model (\ref{equ1.6}) apart of the regression
function $f$
chosen by nature and also knows the dealt class (\ref{equ1.8}).
As an example, we may say that (\ref{equ1.2}) is the lower bound for
the minimax MISE where the
supremum is taken over all regression functions from
${\cal S}(\alpha,Q)$,
and the infimum is taken over all possible dealer-estimators.
% As an example, (\ref{equ1.2}) is the lower bound for the minimax MISE
%of all possible dealer-estimators.
%that know the Sobolev class
%${\cal S}(\alpha,Q)$ and its parameters $\alpha$ and $Q$.
\textit{Oracle-estimator} knows everything that a dealer-estimator does
plus a regression function $f$ chosen by nature.
As we shall see shortly, they may be useful in
suggesting a good estimator.

The context of the paper is as follows. Section~\ref{sec2} presents
main theoretical results.
Section~\ref{sec3} presents the methodology, estimators and a
discussion of
assumptions and results, for a ladder of regression
models where each model is of interest on its own.
Section~\ref{sec4} is devoted to a numerical study.
Proofs, notes and more discussion can be found
in the online Appendix (see~\cite{Efr}).

%s2 #&#
\section{Main results}\label{sec2}

We begin with lower bounds and then show that they are sharp (attainable)
by estimators.

%s2.1 #&#
\subsection{Lower bounds for dealer-estimators}\label{sec2.1}

Using terminology of the \hyperref[sec1]{Intro-} \hyperref[sec1]{duction}, our aim is to
propose a lower
minimax bound for all possible dealer-estimators that know: (i) A
sample of size $n$; (ii) Model (\ref{equ1.6}) where nuisance functions
$g({\mathbf
z})$, $\sigma(x,{\mathbf z})$ and joint design density $p(x,{\mathbf
z})$ are
given and $\varepsilon$ is an independent standard normal random variable;
(iii) Pivot $f_0$, constants $\alpha$ and $Q$ and sequences $\rho_n$
and $M_n$ used to define a family (\ref{equ1.8}). In other words, a
dealer-estimator~$\tilde f^*$ knows everything apart of a regression
function $f$ and
%
%e2.1 #&#
%
\begin{equation}\label{equ2.1}\qquad
\tilde f^*(x):= \tilde f^* \bigl(x,(X,{\mathbf Z},Y)^n,
f_0(x), g({\mathbf z}),p(x,{\mathbf z}),\sigma(x,{\mathbf z}),
\rho_n,M_n,\alpha,Q \bigr).
\end{equation}
Here $(X,{\mathbf Z},Y)^n:= ((X_1,{\mathbf Z}_1,Y_1),\ldots,(X_n,{\mathbf
Z}_n,Y_n))$ denotes a sample.

Please note that, for the dealer who knows the additive component
$g({\mathbf z})$, model (\ref{equ1.6}) is equivalent to the pivotal
model (\ref{equ1.1}).

%
%as2.1 #&#
%
\begin{assumption}\label{assump2.1}
In models (\ref{equ1.1}) and (\ref{equ1.6}) the regression error
$\varepsilon$ is standard normal and independent of $(X,{\mathbf Z})$.
\end{assumption}

%
%as2.2 #&#
%
\begin{assumption}\label{assump2.2}
The joint design\vspace*{1pt} density $p(x,{\mathbf z})$ of
$(X,{\mathbf Z})$ is supported on $[0,1]^{D+1}$, and
$\max(|{\ln(p(x,{\mathbf z}))}|,|{\ln(\sigma (x,{\mathbf z}))}|)$ is
bounded\vspace*{1pt} on $[0,1]^{D+1}$. Function ${\cal I}(x):= \int_{[0,1]^D}
p(x,{\mathbf z}) \sigma^{-2}(x,{\mathbf z})\,d{\mathbf z}$ is Riemann
integrable on $[0,1]$.
\end{assumption}

%
%th2.1 #&#
%
\begin{theorem}\label{theo2.1}
Let Assumptions~\ref{assump2.1} and~\ref{assump2.2} hold. Then
for models (\ref{equ1.1}) and (\ref{equ1.6}) the following lower
minimax bound
for dealer-estimators (\ref{equ2.1}) holds:
%
%e2.2 #&#
%
\begin{eqnarray}\label{equ2.2}\qquad
&&
\inf_{\tilde f^*} \sup_{f \in{\cal F}(f_0,\rho_n,M_n,\alpha,Q)} E\biggl
\{ \int
_0^1 \bigl[\tilde f^*(x) - f(x)
\bigr]^2 \,dx\biggr\}
\nonumber\\[-8pt]\\[-8pt]
&&
\qquad\geq P(\alpha,Q) \biggl[\int_0^1
\frac{dx} {\int_{[0,1]^{D}}p(x,{\mathbf z})
\sigma^{-2}(x,{\mathbf z}) \,d{\mathbf z}} n^{-1}\biggr]^{2\alpha
/(2\alpha+1)} \bigl(1+o_n(1)
\bigr),\nonumber
\end{eqnarray}
where the infimum is taken over all possible
dealer-estimators $\tilde f^*$, and $P(\alpha,Q)$ is defined in (\ref{equ1.3}).
\end{theorem}

Remember that ${\cal S}(\alpha, Q) = {\cal F}(0,\infty,0,\alpha,Q)$,
and this implies that the lower bound also holds for classical
Sobolev classes.
Let us also note that for the case $\sigma(x,{\mathbf z})=\sigma(x)$, with
positive and having
bounded derivatives on $[0,1]$ functions $p(x)$ and $\sigma(x)$
and Sobolev regression functions, the lower bound (\ref{equ2.2}) is
known from~\cite{7}
where it is established via the equivalence (between regression and
filtering in
white noise) principle. In this paper
a different technique of finding a lower bound is employed which
allows us to relax the
assumptions.

The lower bound (\ref{equ2.2}) is challenging for an estimator
to match because the dealer knows everything
apart from an underlying regression function.
Nonetheless, as we shall see shortly,
it is possible to propose an estimator that matches
performance of the best dealer-estimator.

Now let us consider two classical discrete nonparametric regression models,
\textit{Bernoulli and Poisson}~\cite{8,14}.
They may be defined as (\ref{equ1.6}) where
now the distribution of $\varepsilon$ depends on $(X,{\mathbf Z})$
and $Y \in\{0,1\}$ in the Bernoulli case and $Y \in\{0,1,\ldots\}$
in the Poisson case.
Another way to describe these models is
as follows: (i) For Bernoulli regression
we observe a sample from $(X,{\mathbf Z},Y)$ where $Y$ is Bernoulli
and $\operatorname{Pr}(Y=1|X,{\mathbf Z}) = f(X)+g({\mathbf Z})$; (ii)
For Poisson regression we observe a sample from $(X,{\mathbf Z},Y)$
where $Y$ is Poisson and
$E\{Y|X,{\mathbf Z}\} = f(X)+g({\mathbf Z})$.
Furthermore, there is an
extra bona fide restriction on
estimated regression functions. For Bernoulli case a regression
function takes on
values between zero and one, and for Poisson case a regression
function is positive. This is the place where using a pivot and
local/shrinking classes becomes handy.

These two regressions are inherently heteroscedastic because for the
Ber\-noulli regression
%
%e2.3 #&#
%
\begin{equation}\label{equ2.3}
\sigma^2(x,{\mathbf z}):=\sigma_{fg}^2(x,{
\mathbf z}):= \bigl[f(x)+g({\mathbf z})\bigr] \bigl[1-f(x)-g({\mathbf
z})\bigr]
\end{equation}
and for the Poisson regression
%
%e2.4 #&#
%
\begin{equation}\label{equ2.4}
\sigma^2(x,{\mathbf z}):=\sigma_{fg}^2(x,{\mathbf
z}):= f(x)+g({\mathbf z}).
\end{equation}
This is another specific of these regressions because the scale
function contains extra information about the estimand (the regression
function). Can this information help and improve the minimax MISE convergence?
As the following result shows, the answer is ``no.''

%
%th2.2 #&#
%
\begin{theorem}\label{theo2.2}
Consider the above-described
Bernoulli and Poisson regressions. Suppose that Assumption \ref
{assump2.2} holds
with correspondingly defined scale functions (\ref{equ2.3}) or (\ref
{equ2.4}), and
in (\ref{equ1.8}) $M_n \to\infty$ as $n \to\infty$.
For all $(x,{\mathbf z}) \in[0,1]^{D+1}$ it is assumed that the pivot $f_0(x)$,
used in (\ref{equ1.8}), satisfies
$0 < C_* \leq f_0(x) + g({\mathbf z})$ and additionally for the Bernoulli
regression $f_0(x) + g({\mathbf z}) \leq C^{*} <1$. Then
for both regressions,
%
%e2.5 #&#
%
\begin{eqnarray}\label{equ2.5}\qquad
&&
\inf_{\tilde f^*} \sup_{f \in{\cal F}(f_0,\rho_n,M_n,\alpha,Q)\cap
{\cal F}^*(g)} E\biggl\{ \int
_0^1 \bigl[\tilde f^*(x) - f(x)
\bigr]^2 \,dx\biggr\}
\nonumber\\[-8pt]\\[-8pt]
&&\qquad\geq P(\alpha,Q) \biggl[\int_0^1
\frac{dx} {\int_{[0,1]^{D}}p(x,{\mathbf z})
\sigma_{f_0g}^{-2}(x,{\mathbf z}) \,d{\mathbf z}} n^{-1}\biggr]^{2\alpha
/(2\alpha+1)}\bigl(1+o_n(1)
\bigr),\nonumber
\end{eqnarray}
where the infimum is taken over all possible dealer-estimators $\tilde
f^*$, ${\cal F}^*(g)$ is a class of all bona fide $f$ and $P(\alpha,Q)$
is defined in (\ref{equ1.3}).
\end{theorem}

As we see, the lower oracle's bounds are the same for the
normal regression with continuous responses and Bernoulli and Poisson
regressions with discrete responses; this can be explained by the fact that
conditional distributions of responses, given covariates,
belong to exponential families~\cite{6,26}.

The following result, whose proof and a
specific dealer-estimator can be found in the Appendix (see~\cite{Efr}), shows that the
lower bounds are sharp.

%
%th2.3 #&#
%
\begin{theorem}\label{theo2.3}
Lower\vspace*{1pt} bounds (\ref{equ2.2}) and (\ref{equ2.5}) are attainable by a
dealer-estimator $\check f^*(x)$, that is,\vspace*{2pt} $\sup_{f \in{\cal
S}(f_0,\rho_n,M_n,\alpha,Q)} E\{ \int_0^1 [\check f^*(x) - f(x)]^2
\,dx\}$ is not greater than the right-hand sides of (\ref{equ2.2}) and
(\ref{equ2.5}) for the Normal and Bernoulli/Poisson regressions
considered in Theorems~\ref{theo2.1} and~\ref{theo2.2}, respectively.
\end{theorem}

%s2.2 #&#
\subsection{Sharpness of lower bounds for estimators}\label{sec2.2}

Our aim is to show that an estimator can match performance of a
dealer-estimator, that is, an estimator can be adaptive (to underlying
function class and nuisance functions in a regression model) and sharp
minimax.

Introduce: a tensor-product cosine basis $\{\psi_{\mathbf s}({\mathbf v}):=
\prod_{t=1}^D \varphi_{s_t}(v_t)$, ${\mathbf s}:= (s_1,\break\ldots, s_D)
\in\{0,1,\ldots\}^D$, ${\mathbf v}:= (v_1,\ldots,v_D) \in[0,1]^D\}$,
$l_\infty$-norm $\|{\mathbf s}\|_\infty:= \max(s_1,\break\ldots, s_D)$, analytic
function class ${\cal A}:= {\cal A}(\beta_0,\ldots,\beta_D,Q_1):=
\{q(x,{\mathbf z})\dvtx q(x,{\mathbf z}):= \sum_{i,{\mathbf s}} \pi
_{i{\mathbf
s}}\varphi_i(x)\psi_{\mathbf s}({\mathbf z}), |\pi_{i{\mathbf s}}| \leq Q_1
[e^{\beta_0 i} + \sum_{k=1}^D e^{\beta_k s_k}]^{-1},
\min(\beta_0,\beta_1,\ldots,\beta_D) >0,\break Q_1 < \infty\}$ and a
$k$-variate Sobolev class ${\cal S}_k:= \{q(x_1,\ldots,x_k)\dvtx
q(x_1,\ldots,x_k) = \sum_{i_1,\ldots,i_k=0}^\infty q_{i_1,\ldots,i_k}
\prod_{s=1}^k\varphi_{i_s}(x_{i_s}), \sum_{i_1,\ldots,i_k=0}^\infty
[1+\sum_{s=1}^k (2\pi i_s)^{2k}]\times q^2_{i_1,\ldots,i_k} \leq\break Q_2
< \infty\}$; see~\cite{8,29,34}. Parameters of the classes are unknown to
the statistician. In what follows $\nu$'s are generic nonnegative
constants that are used as powers, and $C$'s are generic positive
constants used as factors.

For convenience of future references, let us introduce an array of
assumptions.

%
%as2.3 #&#
%
\begin{assumption}\label{assump2.3}
The following assumptions may be used in
different propositions:

%(a) Scale function $\sigma(x,{\mathbf z})$ is integrable on
%$[0,1]^{D+1}$ and
%$\max(|\ln(p(x,{\mathbf z}))|,$ $|\ln(\sigma^2(x,{\mathbf z}))|) < C$
%for $(x,{\mathbf z}) \in[0,1]^{D+1}$,
(a) Assumption~\ref{assump2.2} holds and regression error $\varepsilon$
satisfies
(\ref{equ1.7}).

(b) Nuisance additive component $g({\mathbf z})$ is bounded and integrable
on $[0,1]^D$ to zero.

%(c) The joint density and the scale function satisfy for $i =1,2$
%p(x,{\mathbf z}) [\sigma^{-2}(x,{\mathbf z})/{\cal I}(x)]^i
%and for any $t >0$
%p(x,{\mathbf z}) [\sigma^{-2}(x,{\mathbf z})/{\cal I}(x)] \varphi_j(x)
(c) The design density satisfies for some $\nu> 0$,
%
%e2.6 #&#
%
\begin{equation}\label{equ2.6}
\sum_{(j,{\mathbf s}) \in\{0,1,\ldots\}^{D+1}} \biggl|\int_{[0,1]^{D+1}}
p(x,{\mathbf
z}) \varphi_j(x) \psi_{\mathbf s}({\mathbf z})\,dx \,d{\mathbf z}\biggr| \leq
C \ln
^\nu(n)
\end{equation}
and for some positive constant $\nu_0$ and any $t > n^{\nu_0}$,
%
%e2.7 #&#
%
\begin{equation}\label{equ2.7}
\sum_{j=0}^\infty\sum
_{\|{\mathbf s}\| > t } \biggl[\int_{[0,1]^{D+1}} p(x,{\mathbf z})
\varphi_j(x) \psi_{\mathbf s}({\mathbf z})\,dx \,d{\mathbf z}
\biggr]^2 \leq C \ln^\nu(n) t^{-D}.
\end{equation}

(d) The $L_2$-approximation of additive component $g({\mathbf z})$
satisfies for any $t > 0$ and some $\nu> 0$
%
%e2.8 #&#
%
\begin{equation}\label{equ2.8}
\sum_{\|{\mathbf s}\| > t } \biggl[\int_{[0,1]^{D}}g({
\mathbf z}) \psi_{\mathbf s}({\mathbf z})\,d{\mathbf z}\biggr]^2 \leq C
t^{-\nu}.
\end{equation}

(e) Two constants,
$c_*$ and $c^*$, are given such that
\mbox{$0 < c_*\,{\leq}\,\sigma^2(x,{\mathbf z}) \,{\leq}\, c^* \,{<}\, \infty$}.\vadjust{\goodbreak}

(f) Design density $p(x,{\mathbf z})$ belongs
to an analytic class ${\cal A}$.

(g) Design density $p(x,{\mathbf z})$ belongs to a $(D+1)$-variate
Sobolev class ${\cal S}_{D+1}$ and nuisance component $g({\mathbf z})$ belongs
to a $D$-variate Sobolev class ${\cal S}_D$.
\end{assumption}

Let us note that: in part (c) a larger class of densities is allowed
for larger~$n$;
if in part (g) we additionally assume that
$g({\mathbf z}) = \sum_{r=1}^D g_r(z_r)$, then
the familiar assumption $g_r \in{\cal S}_1$, $r = 1,\ldots,D$, is
sufficient and the
corresponding proof can be found in the Appendix (see~\cite{Efr}).

The following proposition presents a ladder of settings, each of
interest on its own, for which sharp-minimax and adaptive estimation is
possible. A discussion of the settings and proposed estimators
will be presented in Section~\ref{sec3}.

%
%th2.4 #&#
%
\begin{theorem}\label{theo2.4}
Consider a general additive regression model (\ref{equ1.6}) with the
regression error that may depend on covariates $(X,{\mathbf Z})$ and
satisfying (\ref{equ1.7}). Then for each of the following five sets of
assumptions there exists an estimator that is sharp minimax and matches
performance of the dealer-estimator outlined in Theorem~\ref{theo2.3}:

(1) Additive component $g({\mathbf z})$, design density and scale are
known and
Assumption~\ref{assump2.3}\textup{(a)} holds.

(2) Design density and scale are known and Assumption
\ref{assump2.3}\textup{(a)--(d)} holds.

(3) Design density is known and Assumption
\ref{assump2.3}\textup{(a)--(e)} holds.

(4) Assumption~\ref{assump2.3}\textup{(a), (b), (d), (e), (f)} holds.

(5) Assumption~\ref{assump2.3}\textup{(a), (b), (e), (g)} holds.
\end{theorem}

This result implies the following proposition.

% Regression model (\ref{equ1.6}), considered in Theorem 2.4,
%includes the case of regression errors depending on covariates.
%Furthermore, using family (\ref{equ1.8}), which includes
%local and shrinking Sobolev classes, allows us to consider
%bona fide regression functions. As a result, many classical
%regressions,
%for instance Bernoulli and Poisson ones,
%become particular examples of model (\ref{equ1.6}), and then we get the
%following proposition.

%
%co2.1 #&#
%
\begin{corollary}\label{corollary2.1}
Consider Bernoulli and Poisson regression models discussed in Theorem
\ref{theo2.2}.
%For $(x,{\mathbf z}) \in[0,1]^{D+1}$ suppose that $0 < C_* < f_0(x) +
%g({
%and additionally for Bernoulli regression $f_0(x) + g({\mathbf z})
Then the assertion of Theorem~\ref{theo2.4} holds, and the same estimators
attain the minimax lower bound of Theorem~\ref{theo2.2}.
\end{corollary}

% Next section presents the methodology of construction of sharp minimax
%and adaptive estimators, as well as a detailed discussion of
%assumptions
%and estimators proposed for settings (1)-(5) in Theorem 2.4.
%Let us also note that the pivotal model (\ref{equ1.1}) is considered
%in part
%(1) of Theorem 2.4.

% Estimators are presented and discussed in the next section.

%s3 #&#
\section{Estimation}\label{sec3}

We begin with an explanation of the \textit{methodology} of sharp-minimax
estimation.
Two technical results
are presented for a general regression model. The former is about a
blockwise-shrinkage
oracle-estimator which is adaptive and sharp-minimax. The latter is about
sufficient conditions for an estimator to mimic the oracle.
These two results shed light on the underlying methodology of constructing
sharp-minimax estimators and are of interest on their own.
Then we are presenting specific estimators for each setting considered
in Theorem~\ref{theo2.4}.

To propose a blockwise-shrinkage oracle-estimator,
let $\{B_k,k=1,2,\ldots\}$ be a partition of nonnegative integers
[frequencies of the cosine basis $\{\varphi_j(x)$, $ j=0,1,\ldots\}$]
into nonoverlapping blocks of cardinality (length) $L_k$ such that
$\max(j\dvtx j \in B_k) < \min(j\dvtx j \in B_{k+1})$. The\vadjust{\goodbreak}
blockwise-shrinkage oracle-estimator is defined as
%
%e3.1 #&#
%
\begin{equation}\label{equ3.1}
\hat f^*(x):= \sum_{k=1}^{K_n}
\mu_k \sum_{j \in B_k} \hat
\theta_j \varphi_j(x),
\end{equation}
where $K_n$ is some positive, nondecreasing and integer-valued
sequence,
%
%e3.2 #&#
%
\begin{equation}\label{equ3.2}
\mu_k:= \frac{\Theta_k}{\Theta_k + d n^{-1}}
\end{equation}
is the oracle's shrinkage coefficient for frequencies from the block $B_k$,
$d:= d(p,\sigma)$ is the coefficient of difficulty (\ref{equ1.5}) that
appears in
the lower bounds (\ref{equ2.2}) and (\ref{equ2.5}),
%
%e3.3 #&#
%
\begin{equation}\label{equ3.3}
\Theta_k:= L_k^{-1} \sum
_{j \in B_k} \theta_j^2
\end{equation}
is the Sobolev functional which defines the average energy of $f(x)$
on frequencies from the block $B_k$. A statistic
$\hat\theta_j$, used in (\ref{equ3.1}), is an appropriate estimator of
the Fourier coefficient $\theta_j = \int_0^1 f(x) \varphi_j(x) \,dx$.
For the purposes of this paper, the oracle should be able
to suggest a statistic whose
mean squared error (MSE) satisfies
%
%e3.4 #&#
%
\begin{equation}\label{equ3.4}
E\bigl\{(\hat\theta_j - \theta_j)^2\bigr\}
\leq dn^{-1}\bigl(1+ o_{n}(1) + o_j(1) +
a_j^2 \ln^\nu(n)\bigr),
\end{equation}
where $d$ is defined in (\ref{equ1.5}),
and its squared bias satisfies
%
%e3.5 #&#
%
\begin{equation}\label{equ3.5}
\bigl[E\{\hat\theta_j\} - \theta_j\bigr]^2
\leq n^{-1}\bigl[o_n(1) + o_j(1)+
a_j^2 \ln^\nu(n)\bigr].
\end{equation}
Here and in what follows $\{a_j^2\}$'s are generic summable sequences
($\sum_{j=0}^\infty a_j^2 < \infty$) and $\nu$'s are generic
nonnegative constants that are used in powers.

The following result explains why it is worthwhile to consider
the oracle-estimator (\ref{equ3.1}).

%
%le3.1 #&#
%
\begin{lemma}\label{lemma3.1}
Suppose that in (\ref{equ3.1}) the sequence $K_n$ is large enough to satisfy
the inequality $\sum_{k=1}^{K_n} L_k > n^{1/(2\alpha+1)} \ln(\ln(n+20))$,
and (\ref{equ3.4})--(\ref{equ3.5}) hold. Then
%
%e3.6 #&#
%
\begin{eqnarray}\label{equ3.6}
&&
\sup_{f \in{\cal F}(f_0,\infty,M_n,\alpha,Q)} E\biggl\{\int_0^1
\bigl(\hat f^*(x) - f(x)\bigr)^2 \,dx\biggr\}
\nonumber\\[-8pt]\\[-8pt]
&&\qquad\leq P(\alpha, Q) (d/n)^{2\alpha/(2\alpha+1)}
\bigl(1+o_n(1)\bigr).\nonumber
\end{eqnarray}
\end{lemma}

Let us make several comments about this result: (i) Lemma \ref
{lemma3.1} does not
refer to or is based on a specific regression model; (ii) It was
explained in the \hyperref[sec1]{Introduction} that ${\cal F}(0,\infty,
0, \alpha,Q) =
{\cal S}(\alpha,Q)$ and thus the presented upper bound holds for
classical Sobolev classes; (iii) Using lower bounds of Section~\ref
{sec2} and
relation ${\cal F}(f_0,\rho_n,M_n,\alpha,Q) \subset{\cal
F}(f_0,\infty,M_n,\alpha,Q)$, we conclude that the oracle-estimator is
adaptive and sharp-minimax.\vadjust{\goodbreak}

Now we are in a position to describe the proposed
methodology of developing a data-driven estimator that mimics
the oracle-estimator and is sharp-minimax.

Let us introduce several new sequences and specific blocks used from
now on.
Set: $b_n:= \lfloor\ln(n+20)\rfloor$ where $\lfloor x \rfloor$ denotes
the largest integer which is at most $x$; $c_n:= \lfloor
\ln(b_n)\rfloor$; $m:= \lfloor n/(7c_n)\rfloor$ and it is assumed
that $n$ is large enough so $m >3$; $L_k:= 1$ for $k=1,2,\ldots,
b_n$ and $L_k:= \lfloor(1+b_n^{-1})^k\rfloor$ for $k > b_n$; $K_n$
is the smallest integer such that $\sum_{k=1}^{K_n} L_k > n^{1/3} c_n$;
$B_k:=\{k-1\}$ for
$k=1,2,\ldots,b_n$ and $B_k:=
\{\sum_{s=1}^{k-1}L_{s},\sum_{s=1}^{k-1}L_{s} +1,\ldots,
\sum_{s=1}^{k}L_{s} -1\}$ for $ b_n < k \leq K_n$.

Let us comment on
the specific choice of blocks. The first $b_n$ blocks have unit
lengths, and this choice is
motivated by good performance for small samples.
%and the corresponding study conducted in~\cite{8}.
Then the length of blocks increases geometrically but in such a way
that $L_{k+1}/L_k = 1+o_n(1)$. This choice is motivated by the
asymptotic analysis together
with a good performance for small samples. Let us note that the number
of considered
blocks, $K_n$, is of order $\ln^2(n)$. The largest length of the
blocks, $L_{K_n}$, is of order
$n^{1/3} [\ln(\ln(n))]/\ln(n)$. The total number of estimated low
frequency Fourier coefficients
is of order $n^{1/3}\ln(\ln(n))$. This choice is explained by the fact
that the sum of
not estimated squared Fourier coefficients is of order $o_n(1)
n^{-2\alpha/(2\alpha+1)}$
whenever $\alpha\geq1$. Another way to look at this choice is as follows.
It is known \mbox{\cite{4,6,9,10,11,12,13,14}} that for Sobolev's functions of
order $\alpha$
at most
$n^{1/(2\alpha+1)}c_n$ first Fourier coefficients should be estimated, and
this defines the choice of $K_n$. Furthermore, if it is additionally
known that $\alpha\geq\alpha_0$,
then the total number can be changed to $n^{1/(2\alpha_0+1)} c_n$.

The following proposition explains how to develop
an estimator that matches performance of the oracle.
%Similarly to Lemma 3.1, no underlying model is assumed, and remember
%that
%$C$s are generic positive constants, $a_j^2$s are generic summable
%sequences, $o_s(1) \to0$ as $s \to\infty$, and $d$ is defined in
%(\ref{equ1.5}).

%
%le3.2 #&#
%
\begin{lemma}\label{lemma3.2}
Suppose that there exist two arrays of statistics $\{\hat\Theta_k, k
=1,\ldots,K_n\}$ and $\{\hat\theta_j, j =0,\ldots,\sum_{k=1}^{K_n}
L_k\}$, and a statistic $\hat d$ such that the two arrays and $\hat d$
are mutually independent, the array $\{\hat\theta_j\}$ satisfies
(\ref{equ3.4})--(\ref{equ3.5}), the array $\{\hat\Theta_k\}$ satisfies
for some positive constants $C_1$ and $\nu_1$
%
%e3.7 #&#
%
\begin{equation}\label{equ3.7}
E\bigl\{(\hat\Theta_k - \Theta_k)^4\bigr\}
\leq C_1 L_k^{-1}b_n^{\nu_1}n^{-2}
\bigl(\Theta_k + n^{-1}\bigr)^2,
\end{equation}
and the statistic $\hat d$ satisfies for some constant $C_2 \geq1$
%
%e3.8 #&#
%
\begin{equation}\label{equ3.8}
E\biggl\{\frac{(\hat d - d)^2}{\hat d}\biggr\} = o_n(1),\qquad \hat d \in
\bigl[(C_2b_n)^{-1/4},(C_2b_n)^{1/4}
\bigr] \qquad\mbox{a.s.}
\end{equation}
Then the blockwise-shrinkage estimator
%
%e3.9 #&#
%
\begin{equation}\label{equ3.9}
\hat f(x):= \sum_{k=1}^{K_n}
\frac{\hat\Theta_k}{\hat\Theta_k + \hat
d n^{-1}} I\bigl(\hat\Theta_k > (b_n
n)^{-1}\bigr) \sum_{j \in B_k} \hat
\theta_j \varphi_j(x),
\end{equation}
which mimics the oracle-estimator (\ref{equ3.1}), inherits the sharp-minimax
property of the oracle-estimator described in Lemma~\ref{lemma3.1}, namely
%
%e3.10 #&#
%
\begin{eqnarray}\label{equ3.10}
&&\sup_{f \in{\cal F}(f_0,\infty,M_n,\alpha,Q)} E\biggl\{\int_0^1
\bigl(\hat f(x) - f(x)\bigr)^2 \,dx\biggr\}
\nonumber\\[-8pt]\\[-8pt]
&&\qquad\leq P(\alpha, Q) (d/n)^{2\alpha/(2\alpha+1)}
\bigl(1+o_n(1)\bigr).\nonumber
\end{eqnarray}
\end{lemma}

Now we are in a position to consider settings (1)--(5) of Theorem
\ref{theo2.4}
in turn,
and propose corresponding
statistics $\{\hat\theta_j,\hat\Theta_k,\hat d\}$ used in the
estimator~(\ref{equ3.9}).

%s3.1 #&#
\subsection{Known additive component, design and scale}\label{sec3.1}

This is the case where model (\ref{equ1.6}) transforms into the pivotal model
(\ref{equ1.1}). Because nuisance additive component $g({\mathbf z})$ is known,
without loss of generality we could assume that $g({\mathbf z}) = 0$ or
replace $Y$ by $Y-g({\mathbf Z})$. However, we do not do this because we
would like to indicate what may be done for the case of unknown $g$.
Our idea is to mimic oracle (\ref{equ3.1}) via application of Lemma \ref
{lemma3.2}. To do
this, we need to suggest estimators for Sobolev functionals $\Theta_k$
and Fourier coefficients $\theta_j$; note that the coefficient of
difficulty $d$, defined in (\ref{equ1.5}), is known. Set
%
%e3.11 #&#
%
\begin{equation}\label{equ3.11}\quad
\hat\theta_j:= \frac{1}{n-2m} \sum
_{l=2m+1}^n \frac{[Y_l - \tilde f_{-j}(X_l)- g({\mathbf Z_l})]\sigma
^{-2}(X_l,{\mathbf
Z}_l) \varphi_j(X_l)} {
{\cal I}(X_l)},
\end{equation}
where
%
%e3.12 #&#
%e3.13 #&#
%
\begin{eqnarray}
\label{equ3.12}
{\cal I}(x)&:=& \int_{[0,1]^{D}} p(x,{\mathbf z}) \sigma^{-2}(x,{
\mathbf z}) \,d{\mathbf z},
\\
\label{equ3.13}
\tilde f_{-j}(x) &=& m^{-1}\sum_{l=1}^m
\sum_{i \in{\cal N}_{-j}} \frac{(Y_l - g({\mathbf Z}_l) \varphi
_j(X_l)}{p(X_l)}\varphi_i(x)
\end{eqnarray}
and ${\cal N}_{-j}:= \{0,1,\ldots, b_n\}\setminus\{j\}$.
Note that $\tilde f_{-j}(x)$ estimates $f_{-j}(x):= f(x) - \theta_j
\varphi_j(x)$.
Further,
%
%e3.14 #&#
%
\begin{equation}\label{equ3.14}\quad
\hat\Theta_k:= \frac{2}{m(m-1)}\sum_{m+1 \leq l_1 < l_2\leq2m}
L_k^{-1} \sum_{j \in B_k}
\frac{Y_{l_1} Y_{l_2}\varphi_j(X_{l_1})\varphi_j(X_{l_2})} {
p(X_{l_1},{\mathbf Z}_{l_1})p(X_{l_2},{\mathbf Z}_{l_2})},
\end{equation}
and note that this is U-statistic and unbiased estimate of $\Theta_k$.
This special form of the estimator $\hat\Theta_k$
(it is different from those used in
\cite{4,7,8,9}) implies existence of the fourth moment of
$\hat\Theta_k$ given existence of the fourth moment of the regression error.
Another remark is that we may use the marginal density of $X$ in place
of the joint
design probability density if $ Y_l -g({\mathbf Z}_l)$ is
used in the numerator of (\ref{equ3.14}) in place of $Y_l$.

Let us comment on the estimator (\ref{equ3.11}) of Fourier coefficients
$\theta_j$.
First,
the statistic $\tilde f_{-j}$ is subtracted from the response to\vadjust{\goodbreak}
decrease the MSE. If the subtraction is skipped
then in (\ref{equ3.4}) we would have a larger factor
$(d + \int_0^1 f^2(x) \,dx)$ in place of the wished $d$. Second,
the estimator uses weights (remember the discussion in the \hyperref
[sec1]{Introduction})
%
%e3.15 #&#
%
\begin{equation}\label{equ3.15}
w_l = \frac{\sigma^{-2}(X_l,{\mathbf Z}_l)}{{\cal I}(X_l)} =\frac{\sigma
^{-2}(X_l,{\mathbf Z}_l)}{p(X_l) E\{\sigma^{-2}(X,{\mathbf
Z})|X=X_l\}}.
\end{equation}
This choice of weights yields the wished properties (\ref
{equ3.4})--(\ref{equ3.5}). Note that
if $\sigma(x,{\mathbf z}) = \sigma(x)$, then
weights (\ref{equ3.15}) do not depend on the scale.\vspace*{-2pt}

%
%pr3.1 #&#
%
\begin{proposition}\label{proposition3.1}
Consider\vspace*{1pt} setting (1) of Theorem~\ref{theo2.4}.
Then the block\-wise-shrinkage regression estimator (\ref{equ3.9})
where $\hat\Theta_k$ is defined in (\ref{equ3.14}) and $\hat\theta_j$
in (\ref{equ3.11}), is adaptive to
$(f_0(x), \rho_n, M_n,\alpha,Q)$ and sharp minimax, that is, its MISE satisfies
(\ref{equ3.10}).\vspace*{-2pt}
\end{proposition}

An interesting outcome of the proposition is that no smoothness of the
pivotal regression
function is required for adaptive sharp-minimax estimation, and that
regression error may depend on covariates and have only the
fourth moment.\vspace*{-2pt}

%
%re3.1 #&#
%
\begin{remark}\label{remark3.1}
In Section~\ref{sec4}, where estimators are tested on small samples, we
will study $D$-estimator which is the above-defined estimator without
splitting data. Similarly, all other proposed estimators, when used for
small samples, do not split data.\vspace*{-2pt}
\end{remark}

%s3.2 #&#
\subsection{Known design and scale}\label{sec3.2}

Here the main complication is an unknown additive nuisance component
$g({\mathbf z})$. To mimic the oracle we need to ``remove'' the nuisance
component from the response, and this is a familiar approach in the
additive models literature. As it is shown in the Appendix (see~\cite{Efr}), this
straightforward approach requires an extra assumption about smoothness
of the scale. Because the main topic of the paper is
heteroscedasticity, it is of interest to assume as little as possible
about the scale function. Furthermore, let us remind the reader that
estimation of the scale function is a complicated statistical problem
on its own because quality of estimation depends on smoothness of the
regression function and the scale function~\cite{3}. As a result, even if
for now the scale function is known, it is desirable to assume as
little as possible about its properties and then later use a simple
estimator of the scale.\looseness=-1

The recommended approach is to replace the known
$\sigma^{-2}(x,{\mathbf z})$ by its Fej\'er approximation of order $b_n$,
%
%e3.16 #&#
%
\begin{eqnarray}\label{equ3.16}\quad
&&\sigma^{-2}_{b_n} (x,{\mathbf z})\nonumber\\[-2pt]
&&\qquad:= b_n^{-1}
\sum_{t=0}^{b_n-1} \sum
_{\|(i,{\mathbf s})\|_\infty\leq t} \biggl[\int_{[0,1]^{D+1}}\sigma^{-2}(u,{
\mathbf v})\varphi_i(u) \psi_{\mathbf
s}({\mathbf v})\,du \,d{\mathbf v}
\varphi_i(x) \psi_{\mathbf s}({\mathbf z}) \biggr] \\
&&\qquad=: \sum
_{\|(i,{\mathbf s})\|_\infty< b_n} \eta_{i{\mathbf s}} \varphi_i(x)
\psi_{\mathbf s}({\mathbf z}).\nonumber
\end{eqnarray}
Here $\eta_{i{\mathbf s}}$ are Fej\'er coefficients (note that they
depend on the order $b_n$).
The Fej\'er approximation has
a unique property of
preserving the range of approximated $\sigma^{-2}(x,{\mathbf z})$; see more
about this nice trigonometric approximation in~\cite{2,8,32,36}.
Note that while using Fej\'er's approximation is important,
the choice of its order (here $b_n$) is flexible. We also
replace known ${\cal I}(x)$ by the corresponding approximation
%
%e3.17 #&#
%e3.18 #&#
%
\begin{eqnarray}\label{equ3.17}
{\cal I}_{b_n}(x)&:=& \int_{[0,1]^{D}} p(x,{\mathbf z})
\sigma_{b_n}^{-2}(x,{\mathbf z}) \,d{\mathbf z}
\\
\label{equ3.18}
&=& \sum_{t=0}^\infty\sum
_{\|(i,{\mathbf
s})\|_\infty< b_n} \pi_{t {\mathbf s}} \eta_{i{\mathbf s}}
\varphi_t(x) \varphi_{i}(x),
\end{eqnarray}
where\vspace*{2pt} $\pi_{t{\mathbf s}}:= \int_{[0,1]^{D+1}} p(x,{\mathbf z})
\varphi_t(x) \psi_{\mathbf s}({\mathbf z}) \,dx \,d{\mathbf z}$ are
Fourier coefficients
of the known design density.

Introduce estimates for $f_{-j}(x)$, $\Theta_k$, $g({\mathbf z})$,
and $\theta_j$ in turn. Write
%
%e3.19 #&#
%
\begin{equation}\label{equ3.19}
\tilde f_{-j}(x):= m^{-1} \sum
_{l=1}^m \sum_{i \in{\cal N}_{-j}}
\frac{Y_l \varphi_i(X_l)}{p(X_l,{\mathbf Z}_l)}\varphi_i(x),
\end{equation}
where ${\cal N}_{-j}$ is the same as in (\ref{equ3.13}),
%L_k^{-1} \sum_{i \in B_k}
%{p(X_{l_1}, {\mathbf Z}_{l_1}) p( X_{l_2},{\mathbf Z}_{l_2})},
%
%e3.20 #&#
%
\begin{equation}\label{equ3.20}
\tilde g({\mathbf z}):= m^{-1} \sum_{l=2m+1}^{3m}
\sum_{{\mathbf r} \in{\cal N}_g} \frac{Y_{l} \psi_{\mathbf r}({\mathbf
Z}_l)}{p(X_l,{\mathbf Z}_l)} \psi_{\mathbf
r}({
\mathbf z})
\end{equation}
is the projection series estimator of $g({\mathbf z})$ with
${\cal N}_g:=\{0,1,\ldots, N_g\}^D \setminus\{0\}^D$ and
$N_g:= \lfloor n^{1/D}/b_n^{2/D}\rfloor$,
and
%
%e3.21 #&#
%
\begin{equation}\label{equ3.21}\qquad\quad
\hat\theta_j:= (n-3m)^{-1} \sum
_{l=3m+1}^n \frac{[Y_l - \tilde f_{-j}(X_l) - \tilde g({\mathbf Z}_l)]
\sigma_{b_n}^{-2}(X_l,{\mathbf Z}_l)\varphi_j(X_l)}{{\cal I}_{b_n}(X_l)}.
\end{equation}

%
%pr3.2 #&#
%
\begin{proposition}\label{proposition3.2}
Consider setting (2) of Theorem~\ref{theo2.4}.
Then the estimator (\ref{equ3.9}), with
$\hat\Theta_k$ defined in (\ref{equ3.14}), $\hat\theta_j$ in (\ref
{equ3.21}) and
$\hat d = d$ defined in (\ref{equ1.5}), is adaptive and sharp
minimax, that is, its MISE satisfies (\ref{equ3.10}).
\end{proposition}

Note that no regularity/smoothness of the scale is assumed (it can be
even discontinuous), but we added a very mild assumption (\ref{equ2.8})
on how
well the nuisance additive component can be approximated by the
trigonometric basis.
For instance, (\ref{equ2.8}) holds if in each variable the function
$g(z_1,\ldots,z_D)$
is piecewise Lipschitz of some positive order (note that
Lipschitz functions of order
$\beta< 1$ are often referred to as H\"older functions)~\cite{8}. The
reason why
the proposed Fej\'er approximation of
$\sigma^{-2}(x,{\mathbf z})$
%(note that $[\sigma_{b_n}^{-2}(x,{\mathbf z})]^{-1/2}$ is a
%corresponding
%approximation of the scale $\sigma(x,{\mathbf z})$)
helps is due to the fact that
it is just a weighted sum of
first $b_n$ Fourier terms of $\sigma^{-2}(x,{\mathbf z})$, that is, the
approximation is an extremely smooth
function.
At the same time, the approximation is sufficient for mimicking the
scale and satisfying (\ref{equ3.4})--(\ref{equ3.5}). While this result
is of
interest on its own, it plays a key role in the case of an unknown
scale because it indicates that a rough estimator of the scale may
be sufficient for a sharp-minimax and adaptive estimation.

%s3.3 #&#
\subsection{Known design}\label{sec3.3}

This is a familiar regression problem which includes, as a particular
case, controlled design regressions~\cite{8,12,14,34}. The main issue now
is an appropriate estimation of the scale.
In the assumption for setting (3) of Theorem~\ref{theo2.4}
we still do not impose any restriction on smoothness of an underlying
scale $\sigma(x,{\mathbf z})$ and have not added a new assumption about
the additive nuisance
component $g({\mathbf z})$. On the other hand, we added Assumption \ref
{assump2.3}(e)
which requires
knowledge of the range of the scale function. If the latter is unknown,
then some
information, on how well the scale can be approximated by the trigonometric
basis, is required; see Remark A.3 in the Appendix (see~\cite{Efr}).

Following Lemma~\ref{lemma3.2} we need to propose an estimate of the coefficient
of difficulty $d$ defined in (\ref{equ1.5}), and, following Section \ref
{sec3.2}, we need
to propose an estimate of $\sigma^{-2}_{b_n}(x,{\mathbf z})$. We begin with
the explanation of how to construct an estimate of $d$.
Remember that, according
to Lemma~\ref{lemma3.2}, an estimator should be independent of all
other statistics.
To estimate the scale function we begin with a truncated
projection estimate of $q(x,{\mathbf z}):= f(x) + g({\mathbf z})$,
%
%e3.22 #&#
%
\begin{eqnarray}\label{equ3.22}\qquad
\tilde q_1(x,{\mathbf z})&:=& \max\Biggl(-b_n,\min
\Biggl(b_n,
m^{-1} \sum_{l=3m+1}^{4m}\sum
_{\|(i,{\mathbf r})\|_\infty< b_n} \frac{Y_l\varphi_i(X_l)\psi_{\mathbf
r}({\mathbf Z}_l)}{p(X_l,{\mathbf
Z}_l)}\nonumber\\[-8pt]\\[-8pt]
&&\hspace*{202.5pt}{}\times\varphi_i(x)
\psi_{\mathbf r}({\mathbf z}) \Biggr) \Biggr),\nonumber
\end{eqnarray}
which is used in the following bona fide projection
estimator of $\sigma^2(x,{\mathbf z})$:
%
%e3.23 #&#
%
\begin{equation}\label{equ3.23}
\tilde\sigma_1^2(x,{\mathbf z}):= \max\biggl(c_*, \min
\biggl(c^*, \sum_{\|(i,{\mathbf r})\|_\infty< b_n} \tilde\sigma
_{1i{\mathbf r}}
\varphi_i(x) \psi_{\mathbf r}({\mathbf z})\biggr)\biggr).
\end{equation}
Here $\tilde\sigma_{1i{\mathbf r}} $ is the estimate of
Fourier coefficients $\sigma_{i{\mathbf r}}$ of $\sigma^{2}(x,{\mathbf z})$,
\[
\sigma_{i{\mathbf r}}:= \int_{[0,1]^{D+1}} \sigma^{2}(x,{
\mathbf z}) \varphi_i(x) \psi_{\mathbf r}({\mathbf z}) \,dx \,d{\mathbf z},
\]
and the proposed estimate, motivated by the method of moments, is
%
%e3.24 #&#
%
\begin{equation}\label{equ3.24}
\tilde\sigma_{1i{\mathbf r}}:= m^{-1} \sum
_{l=4m+1}^{5m} \frac{(Y_l -\tilde q_1(X_l,{\mathbf
Z}_l))^2}{p(X_l,{\mathbf Z}_l)}\varphi_i(X_l)
\psi_{\mathbf r}({\mathbf Z}_l).
\end{equation}

With the bona fide estimate (\ref{equ3.23}) of $\sigma^{2}(x,{\mathbf
z})$ at hand,
we plug it in (\ref{equ1.5}) and get
%
%e3.25 #&#
%
\begin{equation}\label{equ3.25}
\tilde d:= \int_0^1 \frac{dx}{\int_{[0,1]^D} p(x,{\mathbf z})
\tilde\sigma_1^{-2}(x,{\mathbf z}) \,d{\mathbf z}}.
\end{equation}

Now we are utilizing the same approach to estimate
$\sigma_{b_n}^{-2}(x,{\mathbf z})$ used by the estimator
$\hat\theta_j$. Remember that, to follow the recipe of Lemma
\ref{lemma3.2}, this estimate should be independent of $\tilde d$. We
define it similarly to (\ref{equ3.22})--(\ref{equ3.24}),
%
%e3.26 #&#
%
\begin{eqnarray}\label{equ3.26}\qquad
\tilde q(x,{\mathbf z})&:=& \max\Biggl(-b_n,\min\Biggl(b_n,
m^{-1} \sum_{l=5m+1}^{6m}\sum
_{\|(i,{\mathbf r})\|_\infty< b_n} \frac{Y_l\varphi_i(X_l)\psi_{\mathbf
r}({\mathbf Z}_l)}{p(X_l,{\mathbf Z}_l)} \nonumber\\[-8pt]\\[-8pt]
&&\hspace*{201.6pt}{}\times\varphi_i(x)
\psi_{\mathbf r}({\mathbf z}) \Biggr) \Biggr)\nonumber
\end{eqnarray}
and
%
%e3.27 #&#
%
\begin{equation}\label{equ3.27}
\tilde\sigma^2(x,{\mathbf z}):= \max\biggl(c_*, \min\biggl(c^*, \sum
_{\|(i,{\mathbf r})\|_\infty< b_n} \tilde\sigma_{i{\mathbf r}}
\varphi_i(x) \psi_{\mathbf r}({\mathbf z})\biggr)\biggr),
\end{equation}
where
%
%e3.28 #&#
%
\begin{equation}\label{equ3.28}
\tilde\sigma_{i{\mathbf r}}:= m^{-1} \sum
_{l=6m+1}^{7m} \frac{(Y_l - \tilde q(X_l, {\mathbf
Z}_l))^2}{p(X_l,{\mathbf Z}_l)}\varphi_i(X_l)
\psi_{\mathbf r}({\mathbf Z}_l).
\end{equation}

Note that now the estimate $\tilde\sigma^2(x,{\mathbf z})$ plays the
role of
$\sigma^2(x,{\mathbf z})$, and then we apply the Fej\'er approximation
(\ref{equ3.16}) to the estimate
(\ref{equ3.27}) and get the estimate of $\sigma_{b_n}^{-2}(x,{\mathbf z})$,
%
%e3.29 #&#
%
\begin{eqnarray}\label{equ3.29}\qquad
&&
\tilde\sigma^{-2}_{b_n} (x,{\mathbf z})\nonumber\\
&&\qquad:=
b_n^{-1} \sum_{t=0}^{b_n-1}
\sum_{\|(i,{\mathbf s})\|_\infty\leq t} \biggl[\int_{[0,1]^{D+1}}\tilde
\sigma^{-2}(u,{\mathbf v})\varphi_i(u) \psi_{\mathbf
s}({
\mathbf v})\,du \,d{\mathbf v}
\varphi_i(x) \psi_{\mathbf s}({\mathbf z}) \biggr] \\
&&\qquad=: \sum
_{\|(i,{\mathbf s})\|_\infty< b_n} \tilde\eta_{i{\mathbf s}} \varphi_i(x)
\psi_{\mathbf s}({\mathbf z}).\nonumber
\end{eqnarray}
Further, following (\ref{equ3.17}) and (\ref{equ3.18}), we define the
plug-in estimate of
${\cal I}_{b_n}(x)$,
%
%e3.30 #&#
%
\begin{eqnarray}\label{equ3.30}
\tilde{\cal I}_{b_n}(x)&:=& \int_{[0,1]^{D}} p(x,{\mathbf z})
\tilde\sigma_{b_n}^{-2}(x,{\mathbf z}) \,d{\mathbf z}
\nonumber\\[-8pt]\\[-8pt]
&=& \sum_{t=0}^\infty\sum
_{\|(i,{\mathbf
s})\|_\infty< b_n} \pi_{t {\mathbf s}} \tilde\eta_{i{\mathbf s}}
\varphi_t(x) \varphi_{i}(x).\nonumber
\end{eqnarray}
%
%where $\pi_{t{\mathbf s}}:= \int_{[0,1]^{D+1}} p(x,{\mathbf z})
%Fourier
%coefficients of the design density.

Finally, mimicking (\ref{equ3.21}), we introduce a new estimator of Fourier
coefficients~$\theta_j$,
%
%e3.31 #&#
%
\begin{equation}\label{equ3.31}\hspace*{25pt}
\hat\theta_j:= (n-7m)^{-1} \sum
_{l=7m+1}^n \frac{[Y_l - \tilde f_{-j}(X_l) - \tilde g({\mathbf Z}_l)]
\tilde\sigma_{b_n}^{-2}(X_l,{\mathbf Z}_l)\varphi_j(X_l)} {
\tilde{\cal I}_{b_n}(X_l)}.
\end{equation}
Here $\tilde f_{-j}$ and $\tilde g$ are estimates (\ref{equ3.19}) and
(\ref{equ3.20}).

%
%pr3.3 #&#
%
\begin{proposition}\label{proposition3.3}
Consider setting (3) of Theorem~\ref{theo2.4}.
Then the estimator (\ref{equ3.9}), with
$\hat\Theta_k$ defined in (\ref{equ3.14}), $\hat\theta_j$ defined in
(\ref{equ3.31}) and
$\hat d = \tilde d$ defined in (\ref{equ3.25}), is adaptive and sharp
minimax, that is, its MISE satisfies (\ref{equ3.10}).
\end{proposition}

Note that a rough estimate of the scale is sufficient, and
no assumption about smoothness of an underlying scale function is made.

%s3.4 #&#
\subsection{Unknown nuisance functions}\label{sec3.4}

Here we relax the last assumption that the design density $p$ is known.
We are considering setting (4) of Theorem~\ref{theo2.4} (with analytic
$p \in
{\cal A}$) and setting (5) (with Sobolev $p \in{\cal S}_{D+1}$)
simultaneously to highlight similarities and differences in proposed
estimators. We will use the indicator $I(p \notin{\cal A})=1$ for
the case of setting (5). Remember that Sobolev classes were discussed
in the \hyperref[sec1]{Introduction}, a nice discussion of analytic
functions can be
found in~\cite{1,25,29,36} and in~\cite{8} they are recommended for modeling and
approximation of a wide variety of densities for the case of small data
sets.

% In the previous subsections it has been explained how to
%deal with unknown additive component and scale if the design is known.
%Furthermore, it has been shown
%that no smoothness of the scale function is required for
%adaptive and sharp-minimax regression estimation. In this subsection
%we will preserve that property.

% The proposed solution is to use the methodology of regression
%estimation developed in the previous subsections.

Because now the design is unknown, all previously defined estimates
become dealer-estimates, and we will use a standard plug-in technique
of using a density estimate in place of an unknown design density. To
follow the recipe of Lemma~\ref{lemma3.2}, we need to plug-in
independent design density estimates in different oracle-estimates, and
this forces us to rewrite one more time all statistics. This is a good
review of what we have done so far. Remember our notation $b_n:=
\lfloor\ln(n+20)\rfloor$, $c_n:= \lfloor\ln(b_n)\rfloor$, and set
${\cal N}_p:=\{0,1,\ldots,N_p\}^{D+1}$, $N_p:= \lfloor b_n c_n\rfloor
I(p \in{\cal A}) +\lfloor n^{1/3(D+1)}\rfloor I(p \notin{\cal A})$.
Note that $N_p$ is a traditional minimax cutoff for the studied
densities. Set $m:= \lfloor n/[(21)c_n]\rfloor$, ${\cal M}_s:= \{(s-1)m
+ 1,(s-1)m + 2,\ldots,sm\}$, and introduce nine identical (but based
on different subsamples) truncated minimax projection density estimates
\cite{5,8}
%
%e3.32 #&#
%
\begin{equation}\label{equ3.32}\qquad
\tilde p_s(x,{\mathbf z}):= \max\biggl(c_n^{-1},m^{-1}
\sum_{l \in{\cal M}_s} \sum_{(i,{\mathbf r}) \in{\cal N}_p}
\varphi_i(X_l)\psi_{{\mathbf r}}({\mathbf
Z}_l)\varphi_i(x) \psi_{\mathbf r}({\mathbf z})\biggr),
\end{equation}
where $s=1,\ldots,9$.
We have
truncated the projection density estimate from below by $c_n^{-1}$
because its reciprocal will be used.

Now we can define statistics used by the proposed estimator.
The first one is the estimator mimicking
dealer-estimator (\ref{equ3.25}) of the coefficient
of difficulty $d$. We begin with mimicking\vadjust{\goodbreak} dealer-estimates (\ref
{equ3.22}) and (\ref{equ3.23})
used in (\ref{equ3.25}). Write
%
%e3.33 #&#
%e3.34 #&#
%e3.35 #&#
%
\begin{eqnarray}\label{equ3.33}
\tilde q_1(x,{\mathbf z})&:=& \max\biggl(-b_n,\min
\biggl(b_n,
m^{-1}\sum_{l \in{\cal M}_{10}}\sum
_{\|(i,{\mathbf r})\|_\infty< b_n} \frac{Y_l\varphi_i(X_l)\psi_{\mathbf
r}({\mathbf Z}_l) }{\tilde p_1(X_l,{\mathbf
Z}_l)}\hspace*{-32pt}\nonumber\\[-8pt]\\[-8pt]
&&\hspace*{196.5pt}{}\times\varphi_i(x)
\psi_{\mathbf r}({\mathbf z}) \biggr) \biggr),\hspace*{-32pt}\nonumber
\\
\label{equ3.34}
\tilde\sigma_{1i{\mathbf r}}&:=& m^{-1} \sum
_{l \in{\cal M}_{11}} \frac{(Y_l -\tilde q_1(X_l,{\mathbf
Z}_l))^2}{\tilde
p_2(X_l,{\mathbf Z}_l)}\varphi_i(X_l)
\psi_{\mathbf r}({\mathbf Z}_l),\hspace*{-32pt}
\\
\label{equ3.35}
&&\hspace*{-33.5pt}\tilde\sigma_1^2(x,{\mathbf z}):= \max\biggl(c_*, \max
\biggl(c^*, \sum_{\|(i,{\mathbf r})\|_\infty< b_n} \tilde\sigma
_{1i{\mathbf r}}
\varphi_i(x) \psi_{\mathbf r}({\mathbf z})\biggr)\biggr).
\end{eqnarray}

These statistics allow us to define the estimate of $d$ [compare with
(\ref{equ3.25})],
%
%e3.36 #&#
%
\begin{equation}\label{equ3.36}
\tilde d:= \int_0^1 \frac{dx}{\int_{[0,1]^D} \tilde p_3(x,{\mathbf z})
\tilde\sigma_1^{-2}(x,{\mathbf z}) \,d{\mathbf z}}.
\end{equation}

Now we consider a number of statistics used to calculate $\hat\theta_j$
and $\hat\Theta_k$. Following (\ref{equ3.19}), set ${\cal N}_{-j}:=
\{\{0,1,\ldots,b_n\}\setminus\{j\}\}I(p \in{\cal A}) +
\{\{0,1,\ldots,\lfloor n^{1/3}\rfloor\}\setminus\{j\}\}I(p \notin{\cal
A})$ and define the estimate of $f_{-j}(x):= f(x) - \theta_j
\varphi_j(x)$ as
%
%e3.37 #&#
%
\begin{equation}\label{equ3.37}
\tilde f_{-j}(x):= m^{-1} \sum
_{l \in{\cal M}_{12}} \sum_{i \in{\cal N}_{-j}}
\frac{Y_l \varphi_i(X_l)}{\tilde p_4(X_l,{\mathbf Z}_l)}\varphi_i(x).
\end{equation}

Following (\ref{equ3.20}), we define the estimate of the additive
nuisance component
$g({\mathbf z})$ as
%
%e3.38 #&#
%
\begin{equation}\label{equ3.38}
\tilde g({\mathbf z}):= m^{-1} \sum_{l \in{\cal M}_{13}}
\sum_{{\mathbf r}
\in{\cal N}_g} \frac{Y_{l} \psi_{\mathbf r}({\mathbf Z}_l)}{\tilde
p_5(X_l,{\mathbf Z}_l)} \psi_{\mathbf r}({
\mathbf z}),
\end{equation}
where
\begin{eqnarray*}
{\cal N}_g&:=&
\bigl\{\bigl\{0,1,\ldots,\bigl\lfloor n^{1/D}/b_n^{2D}\bigr\rfloor\bigr\}^D
\setminus\{0\}^D\bigr\}
I(p \in{\cal A})\\
&&{}+ \bigl\{\bigl\{0,1,\ldots,
\bigl\lfloor n^{1/(3D)}\bigr\rfloor\bigr\}^D \setminus\{0\}^D\bigr\}
I(p \notin{\cal A}).
\end{eqnarray*}

Now we are following (\ref{equ3.26})--(\ref{equ3.30}) and
estimates $\sigma^{-2}_{b_n}(x,{\mathbf z})$ and ${\cal I}_{b_n}(x)$.
Write
%
%e3.39 #&#
%
\begin{eqnarray}\label{equ3.39}\qquad
\tilde q(x,{\mathbf z})&:=& \max\biggl(-b_n,\min\biggl(b_n,
m^{-1} \sum_{l \in{\cal M}_{14}}\sum
_{\|(i,{\mathbf r})\|_\infty< b_n} \frac{Y_l\varphi_i(X_l)\psi_{\mathbf
r}({\mathbf Z}_l)}{\tilde p_6(X_l,
{\mathbf Z}_l)}\nonumber\\[-8pt]\\[-8pt]
&&\hspace*{197pt}{}\times\varphi_i(x)
\psi_{\mathbf r}({\mathbf z}) \biggr) \biggr)\nonumber
\end{eqnarray}
for the estimate of $q(x,{\mathbf z}):= f(x) +g({\mathbf z})$.
This allows us to estimate Fourier coefficients $\sigma_{i{\mathbf r}}$
of the squared scale function by
%
%e3.40 #&#
%
\begin{equation}\label{equ3.40}
\tilde\sigma_{i{\mathbf r}}:= m^{-1} \sum
_{l \in{\cal M}_{15}} \frac{(Y_l - \tilde q(X_l, {\mathbf
Z}_l))^2}{\tilde p_{7}(X_l,{\mathbf Z}_l)} \varphi_i(X_l)
\psi_{\mathbf r}({\mathbf Z}_l).
\end{equation}
Then, following (\ref{equ3.27}), we can define a truncated projection
estimate of the squared scale function
%
%e3.41 #&#
%
\begin{equation}\label{equ3.41}
\tilde\sigma^2(x,{\mathbf z}):= \max\biggl(c_*, \max\biggl(c^*, \sum
_{\|(i,{\mathbf r})\|_\infty< b_n} \tilde\sigma_{i{\mathbf r}}
\varphi_i(x) \psi_{\mathbf r}({\mathbf z})\biggr)\biggr).
\end{equation}

In addition to density estimates (\ref{equ3.32}), let us
introduce three identical (but based on different subsamples)
density estimates.
Set $N_p^*:= N_pI(p \in{\cal A})+
\lfloor n^{1/(2(D+2))} \rfloor I(p \notin{\cal A})$,
${\cal N}^*_p:=\{0,1,\ldots,N_p^*\}^{D+1}$, and for $s=1,2,3$, define
%
%e3.42 #&#
%
\begin{equation}\label{equ3.42}\hspace*{28pt}
\check p_s(x,{\mathbf z}):= \max\biggl(c_n^{-1},m^{-1}
\sum_{l \in{\cal M}_{15+s}} \sum_{(i,{\mathbf r}) \in{\cal N}^*_p}
\varphi_i(X_l)\psi_{{\mathbf r}}({\mathbf
Z}_l)\varphi_i(x) \psi_{\mathbf r}({\mathbf z})\biggr).
\end{equation}
Note that, with respect to (\ref{equ3.32}), the estimate (\ref
{equ3.42}) is changed only for
the case of Sobolev design densities where a larger cutoff
(implying a smaller bias) is used;
a discussion of why the change is needed and what are the other options
can be found in the Appendix (see~\cite{Efr}).

Now we can introduce
estimates for $\sigma^{-2}_{b_n}(x,{\mathbf z})$ and ${\cal
I}_{b_n}(x)$. Following the methodology of (\ref{equ3.29}) and (\ref
{equ3.30}) we set
%
%e3.43 #&#
%
\begin{eqnarray}\label{equ3.43}\qquad
&&
\tilde\sigma^{-2}_{b_n} (x,{\mathbf z})\nonumber\\
&&\qquad:=
b_n^{-1} \sum_{t=0}^{b_n-1}
\sum_{\|(i,{\mathbf s})\|_\infty\leq t} \biggl[\int_{[0,1]^{D+1}}\tilde
\sigma^{-2}(u,{\mathbf v})\varphi_i(u) \psi_{\mathbf
s}({
\mathbf v})\,du \,d{\mathbf v}
\varphi_i(x) \psi_{\mathbf s}({\mathbf z}) \biggr] \\
&&\qquad=: \sum
_{\|(i,{\mathbf
s})\|_\infty< b_n} \tilde\eta_{i{\mathbf s}} \varphi_i(x)
\psi_{\mathbf
s}({\mathbf z})\nonumber
\end{eqnarray}
and [note that the estimate (\ref{equ3.42}) is used]
%
%e3.44 #&#
%
\begin{eqnarray}\label{equ3.44}
\tilde{\cal I}_{b_n}(x)&:=& \int_{[0,1]^{D}} \check
p_{1}(x,{\mathbf z}) \tilde\sigma_{b_n}^{-2}(x,{\mathbf
z}) \,d{\mathbf z}
\nonumber\\[-8pt]\\[-8pt]
&=& \sum_{t=0}^{N_p^*} \sum
_{\|(i,{\mathbf
s})\|_\infty< b_n} \check\pi_{t {\mathbf s}} \tilde\eta_{i{\mathbf s}}
\varphi_t(x) \varphi_{i}(x),\nonumber
\end{eqnarray}
where $\check\pi_{t{\mathbf s}}:= \int_{[0,1]^{D+1}} \check
p_{1}(x,{\mathbf z})
\varphi_t(x) \psi_{\mathbf s}({\mathbf z}) \,dx \,d{\mathbf z}$ are
Fourier coefficients
of the density estimate.

Only for the case of a Sobolev design density do we need to calculate
statistics
%
%e3.45 #&#
%
\begin{eqnarray}\label{equ3.45}\quad
\hat q_{-j,s}(x)&:=& m^{-1} \sum_{l \in{\cal M}_{18+s}}
\biggl[\sum_{i \in{\cal N}_{-j}} \frac{Y_l \varphi_i(X_l)} {
\tilde p_{7+s}(X_l,{\mathbf Z}_l)}\varphi_i(x)
\nonumber\\[-8pt]\\[-8pt]
&&\hspace*{58.5pt}{}+\sum_{{\mathbf r} \in{\cal N}_g}\frac{Y_{l} \psi_{\mathbf
r}({\mathbf
Z}_l)}{\tilde
p_{7+s}(X_l,{\mathbf Z}_l)} \psi_{\mathbf r}({
\mathbf z}) \biggr],\qquad s=1,2.\nonumber
\end{eqnarray}
Here ${\cal N}_{-j}$ and ${\cal N}_g$ are defined above line (\ref
{equ3.37}) and
below line (\ref{equ3.38}), respectively.

This finishes all preliminary calculations. Now we can define a new
estimator for Sobolev functionals,
%
%e3.46 #&#
%
\begin{eqnarray}\label{equ3.46}
\hat\Theta_k&:=& \frac{2}{m(m-1)}
\nonumber\\
&&{}\times\sum
_{l_1,l_2 \in{\cal M}_{21}, l_1
< l_2} L_k^{-1}\sum_{j \in B_k} \frac{ [Y_{l_1} - I(p \notin{\cal A})\hat
q_{-j,1}(X_{l_1},{\mathbf Z}_{l_1})]
\varphi_j(X_{l_1}) } {\check p_2(X_{l_1}, {\mathbf Z}_{l_1}) }
\\
&&\hspace*{113pt}{}\times\frac{ [Y_{l_2}- I(p \notin{\cal A})\hat
q_{-j,2}(X_{l_2},{\mathbf Z}_{l_2})]
\varphi_{j}(X_{l_2})} {
\check p_3( X_{l_2},{\mathbf Z}_{l_2})}\nonumber
\end{eqnarray}
and, mimicking dealer-estimate (\ref{equ3.31}) of
Fourier coefficients $\theta_j$, define
%
%e3.47 #&#
%
\begin{equation}\label{equ3.47}
\hat\theta_j:= (n-21m)^{-1} \sum
_{l=21m+1}^n \frac{[Y_l - \tilde f_{-j}(X_l) - \tilde g({\mathbf Z}_l)]
\tilde\sigma_{b_n}^{-2}(X_l,{\mathbf Z}_l)\varphi_j(X_l)} {
\tilde{\cal I}_{b_n}(X_l)}.\hspace*{-32pt}
\end{equation}
Here $\tilde f_{-j}$, $ \tilde g$, $\tilde\sigma_{b_n}^{-2}$ and
$\tilde{\cal I}_{b_n}$ are defined
in (\ref{equ3.37}), (\ref{equ3.38}), (\ref{equ3.43}) and (\ref
{equ3.44}), respectively.

%
%pr3.4 #&#
%
\begin{proposition}\label{proposition3.4}
Consider settings (4) and (5) of Theorem~\ref{theo2.4}. Assume that
$I(p \in{\cal A})=1$ and $I(p \notin{\cal A})=1$ indicate that settings
(4) and (5) are considered, respectively. Then estimator
(\ref{equ3.9}), with $\hat\Theta_k$ defined in (\ref{equ3.46}),
$\hat\theta_j$ defined in (\ref{equ3.47}) and $\hat d = \tilde d$
defined in (\ref{equ3.36}), is adaptive and sharp minimax, that is, its
MISE satisfies (\ref{equ3.10}).
\end{proposition}

%
%re3.2 #&#
%
\begin{remark}\label{remark3.2} In what follows the proposed
data-driven estimator,
calculated without splitting data and with $I(p \in{\cal A}) = 1$,
is referred to as $S$-estimator.
\end{remark}

Propositions~\ref{proposition3.1} and~\ref{proposition3.4} imply that
the pivotal
model (\ref{equ1.1}) is a fair benchmark for the general additive model
(\ref{equ1.6}), and
this proves the conjecture made in the \hyperref[sec1]{Introduction}.
More discussion, notes and remarks can be found in the Appendix (see~\cite{Efr}).

%s4 #&#
\section{Numerical study}\label{sec4}

We begin with the following Monte Carlo study. The underlying model is
(\ref{equ1.6}) where $D=1$, $g(z) = 0$, the joint design density
$p(x,z) =
I((x,z) \in[0,1]^2)$, the scale function is $\sigma(x,z) = e^{\lambda
z/2}$ and the regression error is standard normal and independent of
the covariates $(X,Z)$. We use $\lambda\in\{1,2,3\}$ and four
sample sizes $n \in\{50,100,200,400\}$. Figure~\ref{fig1} illustrates a
particular simulation with $n=100$ and $\lambda=2$.

We are exploring 4 different estimation procedures with the first two
being sharp-minimax for model (\ref{equ1.6}) and the last two being
sharp-minimax for the model (\ref{equ1.1}) with
$\sigma(x,z)=\sigma(x)$. The first one is $D$-estimator defined in Remark
\ref{remark3.1}. It knows a sample of size $n$ from $(X,Z,Y)$ and all
nuisance functions in the underlying model (\ref{equ1.6}). This
dealer-estimator serves as a benchmark for the data-driven $S$-estimator
defined in Remark~\ref{remark3.2}. The third estimator is the
$E$-estimator of~\cite{7,8} and it was discussed in the
\hyperref[sec1]{Introduction}. $E$-estimator ignores the
heteroscedasticity but nonetheless for the considered experiment with
$g(z)=0$ it is rate-minimax. In what follows an $E$-estimator based on
a sample of size $n$ from $(X,Y)$ will be referred to as the
$En$-estimator. The last estimator is also an $E$-estimator which is
based on a larger sample of size $m$. Namely, the larger sample
includes the sample of size $n$ from $(X,Y)$, available to the three
previous estimators, and then we add extra $m-n$ observations from
$(X,Y)$. Here $m$ is the rounded up $n d_2/d = n \int_0^1 e^{\lambda
z}\,dz \int_0^1 e^{-\lambda z} \,dz$; remember the discussion below
line (\ref{equ1.5}). We will refer to this estimator as the
$Em$-estimator to stress that it is based on a larger sample of size
$m$. The underlying idea of exploring $Em$-estimator is as follows.
According to the asymptotic theory, $D$- and $S$-estimators, based on a
sample of a sufficiently large size $n$, should have the same MISE as
$Em$-estimator which ignores the heteroscedasticity but can use extra
$m-n$ observations. We will test this asymptotic conclusion shortly.

\begin{figure}

\includegraphics{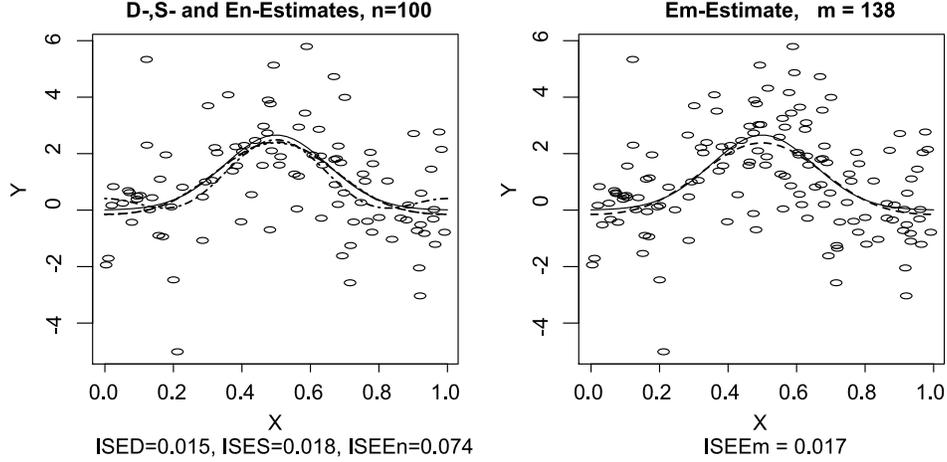}

\caption{Simulated data\vspace*{1pt} according to model (\protect\ref{equ1.6}) with
$f(x)$ being the Normal~\cite{8}, page 18, and shown by the solid line,
$D=1$, $g(z)=0$, $\sigma(x,z) = e^{z}$ and $p(x,z) = I((x,z)
\in[0,1]^2)$. The left scattergram is the same as in the left diagram
of Figure \protect\ref{fig1}, the right scattergram exhibits the same
100 observations
plus 38 additional ones, so the total sample size is $m=138$. All
estimators know that the underlying model is (\protect\ref{equ1.6}),
but only the $D$-estimator knows everything else, except for the
regression function. The left scattergram is overlaid by the
$D$-estimate, $S$-estimate and $En$-estimate shown by the dashed, dotted and
dashed-dotted lines, respectively. The dashed line in the right diagram
shows the $Em$-estimate.}\label{fig2}
\end{figure}

Figure~\ref{fig2} shows us a particular simulation, underlying
regression (the
solid line)
and four estimates (explained in the caption) with their ISE.
For the data, shown in the left diagram, all three estimates
do a very good job under the difficult circumstances, but their ISEs
(denoted as ISED, ISES
and ISEEn, resp.)
reveal that the $D$-estimate
is better than the $S$-estimate, and the $En$-estimate lags behind.
All three estimates give us a fair visualization
of the bell-type and symmetric about 0.5 underlying regression
function. Furthermore, it is practically impossible to see a difference between
the $D$- and $S$-estimates. This highlights the sensitivity of the ISE criterion.
The main issue with the $En$-estimate is in its tails, but they do reflect
the underlying pattern of the shown scattergram (remember that $En$-estimator
knows only the $XY$-scattergram and has no access to observations of $Z$).
The right diagram shows us a scattergram with 38 observations added
from $(X,Y)$. The $Em$-estimate (remember that the same $E$-estimator
is used in the left and right diagrams) yields a much better fit than
the $En$-estimate,
and its ISE (denoted as ISEEm) is close to the ISED and ISES.

%
%f2 #&#
%

For each of 12 particular experiments,
defined by the scale function and the sample size, we conduct 1000
simulations and then calculate average ISE (AISE) for the four estimates.
Table~\ref{table1} presents ratios $R_1:= \mathrm{AISES}/\mathrm
{AISED}$, $R_2:=
\mathrm{AISEEn}/\mathrm{AISES}$ and
$R_3:= \mathrm{AISEEm}/\mathrm{AISED}$.

%t1
%
%t1 #&#
%
\begin{table}
\caption{Results of Monte Carlo simulations}\label{table1}
\begin{tabular*}{\tablewidth}{@{\extracolsep{\fill}}lccccc@{}}
\hline
& & \multicolumn{4}{c@{}}{$\bolds{n}$} \\[-4pt]
& & \multicolumn{4}{c@{}}{\hrulefill} \\
$\bolds{\lambda}$ & & \multicolumn{1}{c}{\textbf{50}}
& \multicolumn{1}{c}{\textbf{100}} & \multicolumn{1}{c}{\textbf{200}}
& \multicolumn{1}{c@{}}{\textbf{400}} \\ \hline
1 & $m$ &54 & 108 & 216 & 432 \\[4pt]
& $R_1,R_2,R_3$ &$1.04,0.87,0.90$ & $1.02,1.08,1.02$
& $1.04,1.12,0.98$ & $1.01,1.21,1.13$ \\
& $R_4,R_5,R_6$ & $1.12,1.14, 1.15$ & $1.08,1.09,1.10$ & $1.07,1.08,1.09$ &
$1.03,1.03,1.04$ \\ [4pt]
2 & $m$ &69 & 138 & 276 & 552 \\[4pt]
& $R_1,R_2,R_3$ &$1.03,0.94,0.79$ & $1.09,1.21,1.01$
& $1.06,1.20,0.95$ &$1.02,1.26,0.96$ \\
& $ R_4,R_5,R_6$ & $1.09,1,11,1,14$ & $1.14,1.15,1.17$ & $1.09,1.10,1.12$ &
$1.04,1.05,1.05$ \\ [4pt]
3 & $m$ &100 & 201 & 403 & 806 \\[4pt]
& $R_1,R_2,R_3$ &$1.61,1.03,0.68$ & $1.11,1.24,0.85$
& $1.09,1.63,0.96$ &$1.06,1.51,0.91$ \\
& $R_4,R_5,R_6$ & $1.78,1.85,1.92$ & $1.18,1.21,1.24$ & $1.12,1.14,1.15$ &
$1.09,1.10,1.11$ \\
\hline
\end{tabular*}
\end{table}

The observed values of ratio $R_1=\mathrm{AISES}/\mathrm{AISED}$
indicate that,
with the exception of
the smallest sample size $n=50$, the proposed data-driven $S$-estimator
does mimic performance
of the
dealer-estimator.
The ratio $R_2 = \mathrm{AISEEn}/\mathrm{AISES}$ shows that even for
the scale
function with
a moderate heteroscedasticity ($\lambda=1$) it may be useful to take
into account the
scale in regression estimation.
Furthermore, the observed values of $R_2$ indicate that a correct usage
of the scale in
regression estimation becomes paramount for regressions with pronounced
heteroscedasticity.
Now let us look at the ratio $R_3= \mathrm{AISEEm}/\mathrm{AISED}$.
The asymptotic theory asserts that the $Em$-estimator,
based on $m$ observations, should have the same MISE as the $D$-estimator
based on $n$ observations (remember Figure~\ref{fig2}).
% This yields that asymptotically $R_3$ should be equal to 1.
As we see, results of the numerical study indicate
that the asymptotic theory sheds light on performance of the
estimators for small samples.
Furthermore, please look at the sample sizes $m$ that make the MISE of
$Em$-estimator
equal to the dealer's MISE. Even for the case $\lambda=1$ we need the
8 percent increase,
and the required sample size doubles for $\lambda=3$.

Now let us repeat simulations three more times using nuisance additive
components $g_1(z) = z-1/2$, $g_2(z) = z^2-1/3$ and $g_3(z) =
z+z^3-3/4$ in place of $g(z) = 0$. We are interested in the effect of a
nuisance component on estimation of~$f$, which can be evaluated via
comparison of performances of the data-driven $S$-estimator and the
$D$-estimator which knows an underlying nuisance component $g_s(z)$.
Results are shown in Table~\ref{table1} via $R_{3+s}:=
\mathrm{AISES}_s/\mathrm{AISED}$, $s=1,2,3$, where $\mathrm{AISES}_s$ is
calculated
for the case of $s$th nuisance component. Note that now $R_1$ serves as
a benchmark for $R_{3+s}$, and we may conclude that $S$-estimator does a
good job in adapting to the presence of a nuisance component.

Overall, the presented numerical results indicate that: (a) Similarly
to~\cite{7,8,27,28},
the asymptotic theory, which takes into account constants, does shed
light on small samples;
(b) It is worthwhile to use the scale in regression estimation whenever
the scale
may depend on auxiliary variables.

\textit{Conclusion}: It is well known that in a nonparametric
heteroscedastic regression the scale function affects the MISE. At the
same time, less is known about optimal use of (or even necessity to
use) the scale function in regression estimation. The pivotal setting,
studied in the paper, is a heteroscedastic regression (\ref{equ1.1})
with a univariate regression function, a multivariate scale and a
normal regression error which is independent of the covariates. For
this setting a sharp-minimax theory of data-driven and adaptive
estimation is developed. The outcome is interesting because, depending
on the scale function, the scale may or may not be recommended for use
by a sharp-minimax regression estimator. Namely, if the scale does not
depend on the auxiliary variable, then a sharp-minimax regression
estimation does not require knowing, using or estimation of the scale,
but otherwise using the scale yields a sharp-minimax MISE. Several
extensions of the pivotal model are also considered: (i) The general
additive model (\ref{equ1.6}) for which model (\ref{equ1.1}) can be
considered as a benchmark. It is shown that the benchmark is fair
meaning that an estimator attains the same minimax MISE for the two
models. Special attention is devoted to assumptions on the nuisance
functions. In particular, it is shown that no smoothness of the scale
is required for the sharp-minimax regression estimation. This is an
important conclusion in light of the known minimax result about the
effect of the smoothness of a regression function on the scale
estimation. Furthermore, the result holds under a mild assumption on
regularity of the multivariate additive component; (ii) The regression
error may not necessarily be normal; it suffices that it has only four
moments, and it may depend on the covariates; (iii) Response may be
discrete with particular examples being classical Bernoulli and Poisson
regressions.
%For these regressions, which are inherently
%heteroscedastic, the sharp-minimax and adaptive theory of regression
%estimation is developed.
A numerical study indicates that the developed sharp-minimax asymptotic theory
sheds light on performance of estimators for small samples.

\section*{Acknowledgments}

The author is grateful for the helpful and constructive comments of the
Editors, Tony Cai and Runze Li, an Associate Editor and two referees.

\begin{supplement}%[id=suppA]
\stitle{Appendix: Notes and proofs}
\slink[doi]{10.1214/13-AOS1126SUPP} %[doi,text={...}] - jei reikia
%suskaldyti doi
\sdatatype{.pdf}
\sfilename{aos1126\_supp.pdf}
\sdescription{Appendix contains proofs and notes.}
\end{supplement}

% imsref loaded by lrinkeviciute, 2013-06-27 10:00:57
% imsref loaded by lrinkeviciute, 2013-06-27 10:03:19

\printaddresses

\end{document}